\documentclass[11pt]{article} 
\usepackage{makeidx}
\usepackage{graphicx}
\usepackage{latexsym}
\usepackage{amsmath}
\usepackage{amscd}
\usepackage{amssymb}
\usepackage{colortab}
\usepackage{array}

\author{St\'ephane Launois 
 \\
\\
{\small{\it Laboratoire de Math\'ematiques - UMR6056, Universit\'e de Reims}}\\
\small{{\it Moulin de la Housse - BP 1039 - 51687 REIMS Cedex 2, France}}\\
\small{e-mail : stephane.launois@univ-reims.fr}}

\title{Rank $t$ $\hc$-primes in quantum matrices.}
\date{ }

\newcommand{\fin}{$\blacksquare$}
\newcommand{\preuve}{\underline{\textit{Proof :}} }
\newcommand{\gc}{ [ \hspace{-0.65mm} [}
\newcommand{\dc}{]  \hspace{-0.65mm} ]}
\newcommand{\ov}{\overline}

\newcommand{\mnk}{O_{q}\left( \mathcal{M}_n(\mathbb{K}) \right) }
\newcommand{\mnc}{O_{q}\left( \mathcal{M}_n(\mathbb{C}) \right)}

\newcommand{\ia}{i,\alpha}

\newcommand{\yia}{Y_{i,\alpha}}

\newcommand{\hc}{\mathcal{H}}

\def\fileversion{97 patch 14}
\def\filedate{1999/12/23}
\csname PSTricksLoaded\endcsname
\PstAtCode{\the\catcode`\@}
\catcode`\@=11\relax
\expandafter\ifx\csname @latexerr\endcsname\relax
\long\def\@ifundefined#1#2#3{\expandafter\ifx\csname
#1\endcsname\relax#2\else#3\fi}
\def\@namedef#1{\expandafter\def\csname #1\endcsname}
\def\@nameuse#1{\csname #1\endcsname}
\def\@eha{%
Your command was ignored.^^J
Type \space I <command> <return> \space to replace
it with another command,^^J
or \space <return> \space to continue without it.}
\def\@spaces{\space\space\space\space}
\def\typeout#1{\immediate\write\@unused{#1}}
\def\@empty{}
\def\@gobble#1{}
\def\@nnil{\@nil}
\def\@ifnextchar#1#2#3{%
\let\@tempe#1\def\@tempa{#2}\def\@tempb{#3}\futurelet\@tempc\@ifnch}
\def\@ifnch{%
\ifx\@tempc\@sptoken
\let\@tempd\@xifnch
\else
\ifx\@tempc\@tempe \let\@tempd\@tempa \else \let\@tempd\@tempb \fi
\fi
\@tempd}
\begingroup
\def\:{\global\let\@sptoken= } \:
\def\:{\@xifnch} \expandafter\gdef\: {\futurelet\@tempc\@ifnch}
\endgroup
\fi
\def\@pstrickserr#1#2{%
\begingroup
\newlinechar`\^^J
\edef\pst@tempc{#2}%
\expandafter\errhelp\expandafter{\pst@tempc}%
\typeout{%
PSTricks error. \space See User's Guide for further information.^^J
\@spaces\@spaces\@spaces\@spaces
Type \space H <return> \space for immediate help.}%
\errmessage{#1}%
\endgroup}
\def\@ehpa{%
Your command was ignored. Default value substituted.^^J
Type \space <return> \space to procede.}
\def\@ehpb{%
Your command was ignored. Will recover best I can.^^J
Type \space <return> \space to procede.}
\def\@ehpc{%
You better fix this before proceding.^^J
See the PSTricks User's Guide or ask your system administrator for help.^^J
Type \space X <return> \space to quit.}
\def\pst@misplaced#1{\@pstrickserr{Misplaced \string#1 command}\@ehpb}
\newdimen\pst@dima
\newdimen\pst@dimb
\newdimen\pst@dimc
\newdimen\pst@dimd
\newdimen\pst@dimg
\newdimen\pst@dimh
\newbox\pst@hbox
\newbox\pst@boxg
\newcount\pst@cnta
\newcount\pst@cntb
\newcount\pst@cntc
\newcount\pst@cntd
\newcount\pst@cntg
\newcount\pst@cnth
\newif\if@pst
\newtoks\pst@toks
\newif\if@star
\def\pst@ifstar#1{%
\@ifnextchar*{\@startrue\def\next*{#1}\next}{\@starfalse#1}}
\def\pst@expandafter#1#2{%
\def\next{#1}%
\edef\@tempa{#2}%
\ifx\@tempa\@empty
\@pstrickserr{Unexpected empty argument!}\@ehpb
\def\@tempa{\@empty}%
\fi
\expandafter\next\@tempa}
\def\pst@dimtonum#1#2{\edef#2{\pst@@dimtonum#1}}
\def\pst@@dimtonum#1{\expandafter\pst@@@dimtonum\the#1}
{\catcode`\p=12 \catcode`\t=12 \global\@namedef{pst@@@dimtonum}#1pt{#1}}
\def\pst@pyth#1#2#3{%
\ifdim#1>#2\pst@@pyth#1#2#3\else\pst@@pyth#2#1#3\fi}
\def\pst@@pyth#1#2#3{%
\ifdim4#1>9#2%
#3=#1\advance#3 .2122#2%
\else
#3=.8384#1\advance#3 .5758#2%
\fi}
\def\pst@divide#1#2#3{%
\pst@@divide{#1}{#2}%
\pst@dimtonum\pst@dimg{#3}}
\def\pst@@divide#1#2{%
\pst@dimg=#1\relax
\pst@dimh=#2\relax
\pst@cntg=\pst@dimh
\pst@cnth=67108863
\pst@@@divide\pst@@@divide\pst@@@divide\pst@@@divide
\divide\pst@dimg\pst@cntg}
\def\pst@@@divide{%
\ifnum
\ifnum\pst@dimg<\z@-\fi\pst@dimg<\pst@cnth
\multiply\pst@dimg\sixt@@n
\else
\divide\pst@cntg\sixt@@n
\fi}
\def\pst@configerr#1{%
\@pstrickserr{\string#1 not defined in pstricks.con}\@ehpc}
\def\pstVerb#1{\pst@configerr\pstVerb}
\def\pstverb#1{\pst@configerr\pstverb}
\def\pstverbscale{\pst@configerr\pstverbscale}
\def\pstrotate{\pst@configerr\pstrotate}
\def\pstheader#1{\pst@configerr\pstheader}
\def\pstdriver{\pst@configerr\pstdriver}
\@ifundefined{pstcustomize}%
{\def\pstcustomize{ \let\pstcustomize\relax}}{}
\input pstricks.con
\newif\ifPSTricks
\PSTrickstrue

\@ifundefined{pst@def}{\def\pst@def#1<#2>{\@namedef{tx@#1}{#2 }}}{}
\@ifundefined{pst@ATH}{\def\pst@ATH<#1>{}}{}
\pstheader{pstricks.pro}
\def\pst@dict{tx@Dict begin }
\def\pst@theheaders{pstricks.pro}
\def\pst@Verb#1{\pstVerb{\pst@dict #1 end}}
\def\tx@Atan{Atan }
\def\tx@Div{Div }
\def\tx@NET{NET }
\def\tx@Pyth{Pyth }
\def\tx@PtoC{PtoC }
\def\tx@PathLength@{PathLength@ }
\def\tx@PathLength{PathLength }
\pst@dimg=\pstunit\relax
\ifdim\pst@dimg=1bp
\def\pst@stp{.996264 dup scale}
\else
\edef\pst@stp{1 \pst@@dimtonum\pst@dimg\space div dup scale}
\fi
\def\tx@STP{STP }
\def\tx@STV{STV }
\def\pst@number#1{\pst@@dimtonum#1\space}
\def\pst@checknum#1#2{%
\edef\next{#1}%
\ifx\next\@empty
\let\pst@num\z@
\else
\expandafter\pst@@checknum\next..\@nil
\fi
\ifnum\pst@num=\z@
\@pstrickserr{Bad number: `#1'. 0 substituted.}\@ehpa
\def#2{0 }%
\else
\edef#2{\ifnum\pst@num=\tw@-\fi\the\pst@cntg.%
\expandafter\@gobble\the\pst@cnth\space}%
\fi}
\def\pst@@checknum{%
\@ifnextchar-%
{\let\pst@num\tw@\expandafter\pst@@@checknum\@gobble}%
{\let\pst@num\@ne\pst@@@checknum}}
\def\pst@@@checknum#1.#2.#3\@nil{%
\afterassignment\pst@@@@checknum\pst@cntg=0#1\relax\@nil
\afterassignment\pst@@@@checknum\pst@cnth=1#2\relax\@nil}
\def\pst@@@@checknum#1\relax\@nil{%
\ifx\@nil#1\@nil\else\let\pst@num\z@\fi}
\def\pst@getnumii#1 #2 #3\@nil{%
\pst@checknum{#1}\pst@tempg
\pst@checknum{#2}\pst@temph}
\def\pst@getnumiii#1 #2 #3 #4\@nil{%
\pst@checknum{#1}\pst@tempg
\pst@checknum{#2}\pst@temph
\pst@checknum{#3}\pst@tempi}
\def\pst@getnumiv#1 #2 #3 #4 #5\@nil{%
\pst@checknum{#1}\pst@tempg
\pst@checknum{#2}\pst@temph
\pst@checknum{#3}\pst@tempi
\pst@checknum{#4}\pst@tempj}
\def\pst@getdimnum#1 #2 #3\@nil{%
\pssetlength\pst@dimg{#1}%
\pst@checknum{#2}\pst@tempg}
\def\pst@getscale#1#2{%
\edef\pst@tempg{#1}%
\ifx\pst@tempg\@none
\def#2{}%
\else
\pst@expandafter\pst@getnumii{#1 #1} {} {} {}\@nil
\ifdim\pst@tempg\p@=\z@
\@pstrickserr{Bad scaling argument `#1'}\@ehpa
\def#2{}%
\else
\ifdim\pst@temph\p@=\z@
\@pstrickserr{Bad scaling argument `#1'}\@ehpa
\def#2{}%
\else
\edef#2{\pst@tempg\space \pst@temph\space scale }%
\fi
\fi
\fi}
\def\pst@getint#1#2{%
\pst@cntg=#1\relax
\edef#2{\the\pst@cntg\space}}
\begingroup
\catcode`\{=12
\catcode`\}=12
\catcode`\[=1
\catcode`\]=2
\gdef\pslbrace[{ ]
\gdef\psrbrace[} ]
\endgroup
\def\@newcolor#1#2{%
\expandafter\edef\csname #1\endcsname{\noexpand\pst@color{#2}}%
\expandafter\edef\csname color@#1\endcsname{#2}%
\ignorespaces}
\def\pst@color#1{%
\def\pst@currentcolor{#1}\pstVerb{#1}\aftergroup\pst@endcolor}
\def\pst@endcolor{\pstVerb{\pst@currentcolor}}
\def\pst@currentcolor{0 setgray}
\def\altcolormode{%
\def\pst@color##1{%
\pstVerb{gsave ##1}\aftergroup\pst@endcolor}%
\def\pst@endcolor{\pstVerb{\pst@grestore}}}
\def\pst@grestore{%
currentpoint
matrix currentmatrix
currentfont
grestore
setfont
setmatrix
moveto}
\def\pst@usecolor#1{\csname color@#1\endcsname\space}
\def\newgray#1#2{%
\pst@checknum{#2}\pst@tempg
\@newcolor{#1}{\pst@tempg setgray}}
\def\newrgbcolor#1#2{%
\pst@expandafter\pst@getnumiii{#2} {} {} {} {}\@nil
\@newcolor{#1}{\pst@tempg \pst@temph \pst@tempi setrgbcolor}}
\def\newhsbcolor#1#2{%
\pst@expandafter\pst@getnumiii{#2} {} {} {} {}\@nil
\@newcolor{#1}{\pst@tempg \pst@temph \pst@tempi sethsbcolor}}
\def\newcmykcolor#1#2{%
\pst@expandafter\pst@getnumiv{#2} {} {} {} {} {}\@nil
\@newcolor{#1}{\pst@tempg \pst@temph \pst@tempi \pst@tempj setcmykcolor}}
\newgray{black}{0}
\newgray{darkgray}{.25}
\newgray{gray}{.5}
\newgray{lightgray}{.75}
\newgray{white}{1}
\newrgbcolor{red}{1 0 0}
\newrgbcolor{green}{0 1 0}
\newrgbcolor{blue}{0 0 1}
\newrgbcolor{yellow}{1 1 0}
\newrgbcolor{cyan}{0 1 1}
\newrgbcolor{magenta}{1 0 1}
\def\psset#1{\@psset#1,\@nil\ignorespaces}
\def\@psset#1,{%
\@@psset#1==\@nil
\@ifnextchar\@nil{\@gobble}{\@psset}}
\def\@@psset#1=#2=#3\@nil{%
\@ifundefined{psset@#1}%
{\@pstrickserr{Graphics parameter `#1' not defined.}\@ehpa}%
{\@nameuse{psset@#1}{#2}}}%
\def\psset@style#1{%
\@ifundefined{pscs@#1}%
{\@pstrickserr{Custom style `#1' undefined}\@ehpa}%
{\@nameuse{pscs@#1}}}
\def\newpsstyle#1#2{\@namedef{pscs@#1}{\psset{#2}}}
\def\@none{none}
\def\pst@getcolor#1#2{%
\@ifundefined{color@#1}%
{\@pstrickserr{Color `#1' not defined}\@eha}%
{\edef#2{#1}}}
\newdimen\psunit \psunit 1cm
\newdimen\psxunit \psxunit 1cm
\newdimen\psyunit \psyunit 1cm
\let\psrunit\psunit
\def\pstunit@off{\let\@psunit\ignorespaces\ignorespaces}
\def\pssetlength#1#2{%
\let\@psunit\psunit
\afterassignment\pstunit@off
#1 #2\@psunit}
\def\psaddtolength#1#2{%
\let\@psunit\psunit
\afterassignment\pstunit@off
\advance#1 #2\@psunit}
\def\pssetxlength#1#2{%
\let\@psunit\psxunit
\afterassignment\pstunit@off
#1 #2\@psunit}
\def\pssetylength#1#2{%
\let\@psunit\psyunit
\afterassignment\pstunit@off
#1 #2\@psunit}
\def\psset@unit#1{%
\pssetlength\psunit{#1}%
\psxunit=\psunit
\psyunit=\psunit}
\def\psset@runit#1{\pssetlength\psrunit{#1}}
\def\psset@xunit#1{\pssetxlength\psxunit{#1}}
\def\psset@yunit#1{\pssetylength\psyunit{#1}}
\def\pst@getlength#1#2{%
\pssetlength\pst@dimg{#1}%
\edef#2{\pst@number\pst@dimg}}
\def\pst@@getlength#1#2{%
\pssetlength\pst@dimg{#1}%
\edef#2{\number\pst@dimg sp}}
\def\pst@getcoor#1#2{\pst@@getcoor{#1}\let#2\pst@coor}
\def\pst@coor{0 0 }
\def\pst@getcoors#1#2{%
\def\pst@aftercoors{\addto@pscode{#1 \pst@coors }#2}%
\def\pst@coors{}%
\pst@@getcoors}
\def\pst@@getcoors(#1){%
\pst@@getcoor{#1}%
\edef\pst@coors{\pst@coor\pst@coors}%
\@ifnextchar({\pst@@getcoors}{\pst@aftercoors}}
\def\pst@getangle#1#2{\pst@@getangle{#1}\let#2\pst@angle}
\def\pst@angle{0 }
\def\cartesian@coor#1,#2,#3\@nil{%
\pssetxlength\pst@dimg{#1}%
\pssetylength\pst@dimh{#2}%
\edef\pst@coor{\pst@number\pst@dimg \pst@number\pst@dimh}}
\def\NormalCoor{%
\def\pst@@getcoor##1{\pst@expandafter\cartesian@coor{##1},\relax,\@nil}%
\def\pst@@getangle##1{%
\pst@checknum{##1}\pst@angle
\edef\pst@angle{\pst@angle \pst@angleunit}}%
\def\psput@##1{\pst@@getcoor{##1}\leavevmode\psput@cartesian}}
\NormalCoor
\def\degrees{\@ifnextchar[{\@degrees}{\def\pst@angleunit{}}}
\def\@degrees[#1]{%
\pst@checknum{#1}\pst@tempg
\edef\pst@angleunit{360 \pst@tempg div mul }%
\ignorespaces}
\def\radians{\def\pst@angleunit{57.2956 mul }}
\def\pst@angleunit{}
\def\SpecialCoor{%
\def\pst@@getcoor##1{%
\begingroup
\pst@activecoor
\xdef\pst@tempg{##1}%
\endgroup
\expandafter\special@coor\pst@tempg||\@nil}%
\def\pst@@getangle##1{%
\begingroup
\pst@activecoor
\xdef\pst@tempg{##1}%
\endgroup
\expandafter\special@angle\pst@tempg\@empty)\@nil}%
\def\psput@##1{\pst@@getcoor{##1}\leavevmode\psput@special}}
\begingroup
\catcode`\|=13
\catcode`\;=13
\catcode`\!=13
\gdef\pst@activecoor{%
\def|{\string|}%
\def;{\string;}%
\def!{\string!}}
\endgroup
\def\special@coor#1|#2|#3\@nil{%
\ifx#3|\relax
\mixed@coor{#1}{#2}%
\else
\special@@coor#1;;\@nil
\fi}
\def\special@@coor#1{%
\ifcat#1a\relax
\def\next{\node@coor#1}%
\else
\ifx#1[\relax
\def\next{\Node@coor[}%
\else
\ifx#1!\relax
\def\next{\raw@coor}%
\else
\def\next{\special@@@coor#1}%
\fi
\fi
\fi
\next}
\def\special@@@coor#1;#2;#3\@nil{%
\ifx#3;\relax
\polar@coor{#1}{#2}%
\else
\cartesian@coor#1,\relax,\@nil
\fi}
\def\mixed@coor#1#2{%
\begingroup
\special@@coor#1;;\@nil
\let\pst@tempa\pst@coor
\special@@coor#2;;\@nil
\xdef\pst@tempg{\pst@tempa pop \pst@coor exch pop }%
\endgroup
\let\pst@coor\pst@tempg}
\def\polar@coor#1#2{%
\pssetlength\pst@dimg{#1}%
\pst@@getangle{#2}%
\edef\pst@coor{\pst@number\pst@dimg \pst@angle \tx@PtoC}}
\def\raw@coor#1;#2\@nil{%
\edef\pst@coor{%
#1 \pst@number\psyunit mul exch \pst@number\psxunit mul exch }}
\def\node@coor#1\@nil{%
\@pstrickserr{You must load `pst-node.tex' to use node coordinates.}\@ehps
\def\pst@coor{0 0 }}
\def\Node@coor{\node@coor}
\def\special@angle#1#2)#3\@nil{%
\ifx!#1\relax
\edef\pst@angle{#2 \pst@angleunit}%
\else
\ifx(#1\relax
\pst@@getcoor{#2}%
\edef\pst@angle{\pst@coor exch \tx@Atan}%
\else
\pst@checknum{#1#2}\pst@angle
\edef\pst@angle{\pst@angle \pst@angleunit}%
\fi
\fi}
\def\Cartesian{%
\def\cartesian@coor##1,##2,##3\@nil{%
\pssetxlength\pst@dimg{##1}%
\pssetylength\pst@dimh{##2}%
\edef\pst@coor{\pst@number\pst@dimg \pst@number\pst@dimh}}%
\@ifnextchar({\Cartesian@}{}}
\def\Cartesian@(#1,#2){%
\pssetxlength\psxunit{#1}%
\pssetylength\psyunit{#2}%
\ignorespaces}
\def\Polar{%
\def\psput@cartesian{\psput@special}%
\def\cartesian@coor##1,##2,##3\@nil{\polar@coor{##1}{##2}}}%
\def\psset@origin#1{%
\pst@@getcoor{#1}%
\edef\psk@origin{\pst@coor \tx@NET }}
\def\psk@origin{}
\newif\ifswapaxes
\def\psset@swapaxes#1{%
\@nameuse{@pst#1}%
\if@pst
\def\psk@swapaxes{-90 rotate -1 1 scale }%
\else
\def\psk@swapaxes{}%
\fi}
\psset@swapaxes{false}
\newif\ifshowpoints
\def\psset@showpoints#1{\@nameuse{showpoints#1}}
\psset@showpoints{false}
\let\pst@setrepeatarrowsflag\relax
\def\psset@border#1{%
\pst@getlength{#1}\psk@border
\pst@setrepeatarrowsflag}
\psset@border{0pt}
\def\psset@bordercolor#1{\pst@getcolor{#1}\psbordercolor}
\psset@bordercolor{white}
\newif\ifpsdoubleline
\def\psset@doubleline#1{%
\@nameuse{psdoubleline#1}%
\pst@setrepeatarrowsflag}
\psset@doubleline{false}
\def\psset@doublesep#1{\def\psdoublesep{#1}}
\psset@doublesep{1.25\pslinewidth}
\def\psset@doublecolor#1{\pst@getcolor{#1}\psdoublecolor}
\psset@doublecolor{white}
\newif\ifpsshadow
\def\psset@shadow#1{%
\@nameuse{psshadow#1}%
\pst@setrepeatarrowsflag}
\psset@shadow{false}
\def\psset@shadowsize#1{\pst@getlength{#1}\psk@shadowsize}
\psset@shadowsize{3pt}
\def\psset@shadowangle#1{\pst@getangle{#1}\psk@shadowangle}
\psset@shadowangle{-45}
\def\psset@shadowcolor#1{\pst@getcolor{#1}\psshadowcolor}
\psset@shadowcolor{darkgray}
\def\pst@repeatarrowsflag{\z@}
\def\pst@setrepeatarrowsflag{%
\edef\pst@repeatarrowsflag{%
\ifdim\psk@border\p@>\z@ 1\else\ifpsdoubleline 1\else
\ifpsshadow 1\else \z@\fi\fi\fi}}
\def\psls@none{}
\newdimen\pslinewidth
\def\psset@linewidth#1{\pssetlength\pslinewidth{#1}}
\psset@linewidth{.8pt}
\def\psset@linecolor#1{\pst@getcolor{#1}\pslinecolor}
\psset@linecolor{black}
\def\psls@solid{0 setlinecap stroke }
\def\pst@missing{%
\z@
\@pstrickserr{Missing number or dimension. 0 substituted}\@ehpa}
\def\psset@dash#1{%
\pst@expandafter\psset@@dash{#1} {\pst@missing} {\pst@missing} {}\@nil
\edef\psk@dash{\pst@number\pst@dimg \pst@number\pst@dimh}}
\def\psset@@dash#1 #2 #3\@nil{%
\pssetlength\pst@dimg{#1}%
\pssetlength\pst@dimh{#2}}
\psset@dash{5pt 3pt}
\newif\ifpsdashadjust
\def\psset@dashadjust#1{\@nameuse{psdashadjust#1}}
\psset@dashadjust{true}
\def\psls@dashed{%
\ifpsdashadjust
\psk@dash \pst@linetype\space \tx@DashLine
\else
[ \psk@dash ] 0 setdash stroke
\fi}
\def\tx@DashLine{DashLine }
\def\psset@dotsep#1{\pst@getlength{#1}\psk@dotsep}
\psset@dotsep{3pt}
\def\psls@dotted{%
\ifpsdashadjust
\psk@dotsep \pst@linetype\space \tx@DotLine
\else
[ 0 \psk@dotsep CLW add ] 0 setdash 1 setlinecap stroke
\fi}
\def\tx@DotLine{DotLine }
\def\psset@linestyle#1{%
\@ifundefined{psls@#1}%
{\@pstrickserr{Line style `#1' not defined}\@eha}%
{\edef\pslinestyle{#1}}}
\psset@linestyle{solid}
\def\psfs@none{}
\def\psset@fillcolor#1{\pst@getcolor{#1}\psfillcolor}
\psset@fillcolor{white}
\def\psfs@solid{\pst@fill{\pst@usecolor\psfillcolor fill}}
\def\psset@hatchwidth#1{\pst@getlength{#1}\psk@hatchwidth}
\psset@hatchwidth{.8pt}
\def\psset@hatchsep#1{\pst@getlength{#1}\psk@hatchsep}
\psset@hatchsep{4pt}
\def\psset@hatchcolor#1{\pst@getcolor{#1}\pshatchcolor}
\psset@hatchcolor{black}
\def\psset@hatchangle#1{\pst@getangle{#1}\psk@hatchangle}
\psset@hatchangle{45}
\def\pst@linefill{%
\psk@hatchangle rotate
\psk@hatchwidth SLW
\pst@usecolor\pshatchcolor
\psk@hatchsep \tx@LineFill}
\def\psfs@vlines{\pst@fill\pst@linefill}
\@namedef{psfs@vlines*}{\psfs@solid \psfs@vlines}
\def\psfs@hlines{\pst@fill{90 rotate \pst@linefill}}
\@namedef{psfs@hlines*}{\psfs@solid \psfs@hlines}
\def\psfs@crosshatch{\psfs@vlines \psfs@hlines}
\@namedef{psfs@crosshatch*}{\psfs@solid \psfs@vlines \psfs@hlines}
\def\tx@LineFill{LineFill }
\def\psset@fillstyle#1{%
\edef\pst@tempg{#1}\def\pst@temph{none}%
\ifx\pst@tempg\pst@temph
\let\psk@fillstyle\relax
\else
\@ifundefined{psfs@#1}%
{\@pstrickserr{Undefined fill style: `#1'}\@eha}%
{\edef\psk@fillstyle{\expandafter\noexpand\csname psfs@#1\endcsname}}%
\fi}
\def\psset@addfillstyle#1{%
\@ifundefined{psfs@#1}%
{\@pstrickserr{Undefined fill style: `#1'}\@eha}%
{\edef\psk@fillstyle{%
\expandafter\noexpand\psk@fillstyle
\expandafter\noexpand\csname psfs@#1\endcsname}}}
\psset@fillstyle{none}
\def\psset@arrows#1{%
\begingroup
\pst@activearrows
\xdef\pst@tempg{#1}%
\endgroup
\expandafter\psset@@arrows\pst@tempg\@empty-\@empty\@nil
\if@pst\else
\@pstrickserr{Bad arrows specification: #1}\@ehpa
\fi}
\def\psset@@arrows#1-#2\@empty#3\@nil{%
\@psttrue
\def\next##1,#1-##2,##3\@nil{\def\pst@tempg{##2}}%
\expandafter\next\pst@arrowtable,#1-#1,\@nil
\@ifundefined{psas@\pst@tempg}%
{\@pstfalse\def\psk@arrowA{}}%
{\let\psk@arrowA\pst@tempg}%
\@ifundefined{psas@#2}%
{\@pstfalse\def\psk@arrowB{}}%
{\def\psk@arrowB{#2}}}
\def\psk@arrowA{}
\def\psk@arrowB{}
\def\pst@arrowtable{,<->,<<->>,>-<,>>-<<,(-),[-]}
\begingroup
\catcode`\<=13
\catcode`\>=13
\catcode`\|=13
\gdef\pst@activearrows{\def<{\string<}\def>{\string>}\def|{\string|}}
\endgroup
\def\tx@BeginArrow{BeginArrow }
\def\tx@EndArrow{EndArrow }
\def\psset@arrowscale#1{\pst@getscale{#1}\psk@arrowscale}
\psset@arrowscale{1}
\def\psset@arrowsize#1{%
\pst@expandafter\pst@getdimnum{#1} 0 {} {}\@nil
\edef\psk@arrowsize{\pst@number\pst@dimg \pst@tempg}}
\psset@arrowsize{1.5pt 2}
\def\psset@arrowlength#1{\pst@checknum{#1}\psk@arrowlength}
\psset@arrowlength{1.4}
\def\psset@arrowinset#1{\pst@checknum{#1}\psk@arrowinset}%
\psset@arrowinset{.4}
\def\tx@Arrow{Arrow }
\@namedef{psas@>}{%
false \psk@arrowinset \psk@arrowlength \psk@arrowsize \tx@Arrow}
\@namedef{psas@>>}{%
false \psk@arrowinset \psk@arrowlength \psk@arrowsize \tx@Arrow
0 h T
gsave
newpath
false \psk@arrowinset \psk@arrowlength \psk@arrowsize \tx@Arrow
CP
grestore
CP newpath moveto
2 copy
L
stroke
moveto}
\@namedef{psas@<}{%
true \psk@arrowinset \psk@arrowlength \psk@arrowsize \tx@Arrow}
\@namedef{psas@<<}{%
true \psk@arrowinset \psk@arrowlength \psk@arrowsize \tx@Arrow
CP newpath moveto 0 a neg L stroke 0 h neg T
false \psk@arrowinset \psk@arrowlength \psk@arrowsize \tx@Arrow}
\def\psset@tbarsize#1{%
\pst@expandafter\pst@getdimnum{#1} 0 {} {}\@nil
\edef\psk@tbarsize{\pst@number\pst@dimg \pst@tempg}}
\psset@tbarsize{2pt 5}
\def\tx@Tbar{Tbar }
\@namedef{psas@|}{\psk@tbarsize \tx@Tbar}
\@namedef{psas@|*}{0 CLW -2 div T \psk@tbarsize \tx@Tbar}
\@namedef{psas@>|}{%
\psk@tbarsize \tx@Tbar
0 CLW 2 div T
newpath
false \psk@arrowinset \psk@arrowlength \psk@arrowsize \tx@Arrow}
\@namedef{psas@>|*}{%
0 CLW -2 div T
\psk@tbarsize \tx@Tbar
0 CLW 2 div T
newpath
false \psk@arrowinset \psk@arrowlength \psk@arrowsize \tx@Arrow}
\edef\pst@arrowtable{\pst@arrowtable,|<*->|*,|<->|}
\def\psset@bracketlength#1{\pst@checknum{#1}\psk@bracketlength}
\psset@bracketlength{.15}
\def\tx@Bracket{Bracket }
\@namedef{psas@]}{\psk@bracketlength \psk@tbarsize \tx@Bracket}
\def\psset@rbracketlength#1{\pst@checknum{#1}\psk@rbracketlength}
\psset@rbracketlength{.15}
\def\tx@RoundBracket{RoundBracket }
\@namedef{psas@)}{\psk@rbracketlength \psk@tbarsize \tx@RoundBracket}
\def\psas@c{1 \psas@@c}
\def\psas@cc{0 CLW 2 div T 1 \psas@@c}
\def\psas@C{2 \psas@@c}
\def\psas@@c{%
setlinecap
0 0 moveto
0 CLW 2 div L
stroke
0 0 moveto}
\def\psas@{}
\psset@arrows{-}
\def\tx@SD{SD }
\def\tx@EndDot{EndDot }
\def\psas@oo{{\pst@usecolor\psfillcolor true} true \psk@dotsize \tx@EndDot}
\def\psas@o{{\pst@usecolor\psfillcolor true} false \psk@dotsize \tx@EndDot}
\@namedef{psas@**}{{false} true \psk@dotsize \tx@EndDot}
\@namedef{psas@*}{{false} false \psk@dotsize \tx@EndDot}
\def\pst@par{}
\def\addto@par#1{%
\ifx\pst@par\@empty
\def\pst@par{#1}%
\else
\expandafter\def\expandafter\pst@par\expandafter{\pst@par,#1}%
\fi}
\def\addbefore@par#1{%
\ifx\pst@par\@empty
\def\pst@par{#1}%
\else
\toks@{#1}%
\pst@toks\expandafter{\pst@par}%
\edef\pst@par{\the\toks@,\the\pst@toks}%
\fi}
\def\use@par{%
\ifx\pst@par\@empty\else
\expandafter\@psset\pst@par,\@nil
\def\pst@par{}%
\fi}
\def\pst@object#1{%
\pst@ifstar{%
\@ifnextchar[%
{\pst@@object{#1}}%
{\def\pst@par{}\@nameuse{#1@i}}}}
\def\pst@@object#1[#2]{%
\def\pst@par{#2}%
\@ifnextchar+{\@nameuse{#1@i}}{\@nameuse{#1@i}}}
\def\newpsobject#1#2#3{%
\@ifundefined{#2@i}%
{\@pstrickserr{Graphics object `#2' not defined}\@eha}{%
\@namedef{#1}{\pst@object{#1}}%
\@namedef{#1@i}{\addbefore@par{#3}\@nameuse{#2@i}}}%
\ignorespaces}
\def\pst@getarrows#1{\@ifnextchar({#1}{\pst@@getarrows{#1}}}
\def\pst@@getarrows#1#2{\addto@par{arrows=#2}#1}
\def\begin@ClosedObj{%
\leavevmode
\pst@killglue
\begingroup
\use@par
\solid@star
\ifpsdoubleline \pst@setdoublesep \fi
\init@pscode}
\def\end@ClosedObj{%
\ifpsshadow \pst@closedshadow \fi
\ifdim\psk@border\p@>\z@ \pst@addborder \fi
\psk@fillstyle
\pst@stroke
\ifpsdoubleline \pst@doublestroke \fi
\ifshowpoints
\pst@OpenShowPoints
\fi
\use@pscode
\endgroup
\ignorespaces}
\def\begin@OpenObj{%
\begin@ClosedObj
\let\pst@linetype\pst@arrowtype
\pst@addarrowdef}
\def\begin@AltOpenObj{%
\begin@ClosedObj
\def\pst@repeatarrowsflag{\z@}%
\def\pst@linetype{0}}
\def\end@OpenObj{%
\ifpsshadow \pst@openshadow \fi
\ifdim\psk@border\p@>\z@ \pst@addborder \fi
\psk@fillstyle
\pst@stroke
\ifpsdoubleline \pst@doublestroke \fi
\ifnum\pst@repeatarrowsflag>\z@ \pst@repeatarrows \fi
\ifshowpoints \pst@OpenShowPoints \fi
\use@pscode
\endgroup
\ignorespaces}
\def\begin@SpecialObj{%
\leavevmode
\pst@killglue
\begingroup
\use@par
\init@pscode}
\def\end@SpecialObj{%
\use@pscode
\endgroup
\ignorespaces}
\def\pst@code{}%
\def\init@pscode{%
\addto@pscode{%
\pst@number\pslinewidth SLW
\pst@usecolor\pslinecolor}}
\def\addto@pscode#1{\xdef\pst@code{\pst@code#1\space}}
\def\use@pscode{%
\pstverb{%
\pst@dict
\tx@STP
\pst@newpath
\psk@origin
\psk@swapaxes
\pst@code
end}%
\gdef\pst@code{}}
\def\pst@newpath{newpath }
\def\pst@@killglue{\unskip\ifdim\lastskip>\z@\expandafter\pst@@killglue\fi}
\def\KillGlue{\let\pst@killglue\pst@@killglue}
\def\DontKillGlue{\let\pst@killglue\relax}
\DontKillGlue
\def\solid@star{%
\if@star
\pslinewidth=\z@
\psdoublelinefalse
\def\pslinestyle{none}%
\def\psk@fillstyle{\psfs@solid}%
\let\psfillcolor\pslinecolor
\fi}
\def\pst@setdoublesep{%
\pst@getlength\psdoublesep\psdoublesep
\pslinewidth=2\pslinewidth
\advance\pslinewidth\psdoublesep\p@
\let\pst@setdoublesep\relax}
\def\tx@Shadow{Shadow }
\def\pst@closedshadow{%
\addto@pscode{%
gsave
\psk@shadowsize \psk@shadowangle \tx@PtoC
\tx@Shadow
\pst@usecolor\psshadowcolor
gsave fill grestore
stroke
grestore
gsave
\pst@usecolor\psfillcolor
gsave fill grestore
stroke
grestore}}
\def\pst@openshadow{%
\addto@pscode{%
gsave
\psk@shadowsize \psk@shadowangle \tx@PtoC
\tx@Shadow
\pst@usecolor\psshadowcolor
\ifx\psk@fillstyle\relax\else
gsave fill grestore
\fi
stroke}%
\pst@repeatarrows
\addto@pscode{grestore}
\ifx\psk@fillstyle\relax\else
\addto@pscode{%
gsave
\pst@usecolor\psfillcolor
gsave fill grestore
stroke
grestore}
\fi}
\def\pst@addborder{%
\addto@pscode{%
gsave
\psk@border 2 mul
CLW add SLW
\pst@usecolor\psbordercolor
stroke
grestore}}
\def\pst@stroke{%
\ifx\pslinestyle\@none\else
\addto@pscode{%
gsave
\pst@number\pslinewidth SLW
\pst@usecolor\pslinecolor
\@nameuse{psls@\pslinestyle}
grestore}%
\fi}
\def\pst@fill#1{\addto@pscode{gsave #1 grestore}}%
\def\pst@doublestroke{%
\addto@pscode{%
gsave
\psdoublesep SLW
\pst@usecolor\psdoublecolor
stroke
grestore}}
\def\pst@arrowtype{%
\ifx\psk@arrowB\@empty 0 \else -2 \fi
\ifx\psk@arrowA\@empty 0 \else -1 \fi
add}
\def\pst@addarrowdef{%
\addto@pscode{%
/ArrowA {
\ifx\psk@arrowA\@empty
\pst@oplineto
\else
\pst@arrowdef{A}
moveto
\fi
} def
/ArrowB {
\ifx\psk@arrowB\@empty \else \pst@arrowdef{B} \fi
} def}}
\def\pst@arrowdef#1{%
\ifnum\pst@repeatarrowsflag>\z@
/Arrow#1c [ 6 2 roll ] cvx def Arrow#1c
\fi
\tx@BeginArrow
\psk@arrowscale
\@nameuse{psas@\@nameuse{psk@arrow#1}}
\tx@EndArrow}
\def\pst@repeatarrows{%
\addto@pscode{%
gsave
\ifx\psk@arrowA\@empty\else
ArrowAc ArrowA pop pop
\fi
\ifx\psk@arrowB\@empty\else
ArrowBc ArrowB pop pop pop pop
\fi
grestore}}
\def\pst@OpenShowPoints{%
\addto@pscode{%
gsave
\psk@dotsize
\@nameuse{psds@\psk@dotstyle}
newpath
Points aload length 2 div 2 sub cvi /N ED
N 0 ge
{ \ifx\psk@arrowA\@empty
Dot
\else
pop pop
\fi
N { Dot } repeat
\ifx\psk@arrowB\@empty
Dot
\else
pop pop
\fi }
{ N 2 mul { pop } repeat }
ifelse
grestore}}
\def\pscustom{\pst@object{pscustom}}
\long\def\pscustom@i#1{%
\begin@SpecialObj
\solid@star
\let\pst@ifcustom\iftrue
\let\begin@ClosedObj\begin@CustomObj
\let\end@ClosedObj\endgroup
\def\begin@OpenObj{\begin@CustomObj\pst@addarrowdef}%
\let\end@OpenObj\endgroup
\let\begin@AltOpenObj\begin@CustomObj
\def\begin@SpecialObj{%
\begingroup
\pst@misplaced{special graphics object}%
\def\addto@pscode####1{}
\let\end@SpecialObj\endgroup}%
\def\@multips(##1)(##2)##3##4{\pst@misplaced\multips}%
\def\psclip##1{\pst@misplaced\psclip}%
\def\pst@repeatarrowsflag{\z@}%
\let\pst@setrepeatarrowsflag\relax
\showpointsfalse
\let\showpointstrue\relax
\def\pst@linetype{\pslinetype}%
\let\psset@liftpen\psset@@liftpen
\psset@liftpen{\z@}%
\def\pst@cp{/currentpoint load stopped pop }%
\def\pst@oplineto{/lineto load stopped { moveto } if }%
\def\pst@optcp##1##2{%
\ifnum##1=\z@\def##2{/currentpoint load stopped { 0 0 } if }\fi}%
\let\caddto@pscode\addto@pscode
\def\cuse@par##1{{\use@par##1}}%
\the\pst@customdefs
\setbox\pst@hbox=\hbox{#1}%
\psk@fillstyle
\pst@stroke
\end@SpecialObj}
\def\begin@CustomObj{%
\begingroup
\use@par
\addto@pscode{%
\pst@number\pslinewidth SLW
\pst@usecolor\pslinecolor}}
\def\pst@oplineto{moveto }
\def\pst@cp{}
\def\pst@optcp#1#2{}
\def\psset@liftpen#1{}
\def\psset@@liftpen#1{%
\ifcase#1\relax
\def\psk@liftpen{\z@}%
\def\pst@cp{/currentpoint load stopped pop }%
\def\pst@oplineto{/lineto load stopped { moveto } if }%
\or
\def\psk@liftpen{1}%
\def\pst@cp{}%
\def\pst@oplineto{/lineto load stopped { moveto } if }%
\or
\def\psk@liftpen{2}%
\def\pst@cp{}%
\def\pst@oplineto{moveto }%
\fi}
\psset@liftpen{0}
\def\psk@liftpen{-1}
\def\psset@linetype#1{%
\pst@getint{#1}\pslinetype
\ifnum\pst@dimg<-3
\@pstrickserr{linetype must be greater than -3}\@ehpa
\def\pslinetype{0}%
\fi}
\psset@linetype{0}
\def\caddto@pscode#1{%
\@pstrickserr{Command can only be used in \string\pscustom}\@ehpa}
\let\cuse@par\caddto@pscode
\def\tx@MSave{%
/msavemtrx
[ tx@Dict /msavemtrx known { msavemtrx aload pop } if CM ]
def }
\def\tx@MRestore{%
tx@Dict /msavemtrx known { length 0 gt } { false } ifelse
{ /msavematrx [ msavematrx aload pop setmatrix ] def }
if }
\newtoks\pst@customdefs
\pst@customdefs{%
\def\newpath{\addto@pscode{newpath}}%
\def\moveto(#1){\pst@@getcoor{#1}\addto@pscode{\pst@coor moveto}}%
\def\closepath{\addto@pscode{closepath}}%
\def\gsave{\begingroup\addto@pscode{gsave}}%
\def\grestore{\endgroup\addto@pscode{grestore}}%
\def\translate(#1){\pst@@getcoor{#1}\addto@pscode{\pst@coor translate}}%
\def\rotate#1{\pst@@getangle{#1}\addto@pscode{\pst@angle rotate}}%
\def\scale#1{\pst@getscale{#1}\pst@tempg\addto@pscode{\pst@tempg}}%
\def\msave{\addto@pscode{\tx@MSave}}%
\def\mrestore{\addto@pscode{\tx@MRestore}}%
\def\swapaxes{\addto@pscode{-90 rotate -1 1 scale}}%
\def\stroke{\pst@object{stroke}}%
\def\fill{\pst@object{fill}}%
\def\openshadow{\pst@object{openshadow}}%
\def\closedshadow{\pst@object{closedshadow}}%
\def\movepath(#1){\pst@@getcoor{#1}\addto@pscode{\pst@coor \tx@Shadow}}%
\def\lineto{\pst@onecoor{lineto}}%
\def\rlineto{\pst@onecoor{rlineto}}%
\def\curveto{\pst@threecoor{curveto}}%
\def\rcurveto{\pst@threecoor{rcurveto}}%
\def\code#1{\addto@pscode{#1}}%
\def\coor(#1){\pst@@getcoor{#1}\addto@pscode\pst@coor\@ifnextchar({\coor}{}}%
\def\rcoor{\pst@getcoors{}{}}%
\def\dim#1{\pssetlength\pst@dimg{#1}\addto@pscode{\pst@number\pst@dimg}}%
\def\setcolor#1{%
\@ifundefined{color@#1}{}{\addto@pscode{\use@color{#1}}}}%
\def\arrows#1{{\psset@arrows{#1}\pst@addarrowdef}}%
\let\file\pst@rawfile
} 
\def\closedshadow@i{\cuse@par\pst@closedshadow}
\def\openshadow@i{\cuse@par\pst@openshadow}
\def\stroke@i{\cuse@par\pst@stroke}%
\def\fill@i{\cuse@par\psk@fillstyle}%
\def\pst@onecoor#1(#2){%
\pst@@getcoor{#2}%
\addto@pscode{\pst@coor #1}}
\def\pst@threecoor#1(#2)#3(#4)#5(#6){%
\begingroup
\pst@getcoor{#2}\pst@tempa
\pst@getcoor{#4}\pst@tempb
\pst@getcoor{#6}\pst@tempc
\addto@pscode{\pst@tempa \pst@tempb \pst@tempc #1}%
\endgroup}
\def\pst@rawfile#1{%
\begingroup
\def\do##1{\catcode`##1=12\relax}"
\dospecials
\catcode`\%=14
\pst@@rawfile{#1}%
\endgroup}
\def\pst@@rawfile#1{%
\immediate\openin1 #1
\ifeof1
\@pstrickserr{File `#1' not found}\@ehpa
\else
\immediate\read1 to \pst@tempg
\loop
\ifeof1 \@pstfalse\else\@psttrue\fi
\if@pst
\addto@pscode\pst@tempg
\immediate\read1 to \pst@tempg
\repeat
\fi
\immediate\closein1\relax}
\def\tx@NArray{NArray }
\def\tx@NArray{NArray }
\def\tx@Line{Line }
\def\tx@Arcto{Arcto }
\def\tx@CheckClosed{CheckClosed }
\def\tx@Polygon{Polygon }
\def\psset@gangle#1{\pst@getangle{#1}\psk@gangle}
\psset@gangle{0}
\def\tx@Diamond{Diamond }
\def\psdiamond{\pst@object{psdiamond}}
\def\psdiamond@i(#1){%
\@ifnextchar({\psdiamond@ii(#1)}{\psdiamond@ii(0,0)(#1)}}
\def\psdiamond@ii(#1)(#2){%
\begin@ClosedObj
\pst@getcoor{#1}\pst@tempa
\pst@getcoor{#2}\pst@tempb
\addto@pscode{%
\psline@iii
pop
\psk@dimen
\pst@tempb
\psk@gangle
\pst@tempa
\tx@Diamond}%
\def\pst@linetype{4}%
\end@ClosedObj}
\def\tx@Triangle{Triangle }
\def\pstriangle{\pst@object{pstriangle}}
\def\pstriangle@i(#1){%
\@ifnextchar({\pstriangle@ii(#1)}{\pstriangle@ii(0,0)(#1)}}
\def\pstriangle@ii(#1)(#2){%
\begin@ClosedObj
\pst@getcoor{#1}\pst@tempa
\pst@getcoor{#2}\pst@tempb
\addto@pscode{%
\psline@iii
pop
\psk@dimen
\pst@tempb
\psk@gangle
\pst@tempa
\tx@Triangle}%
\def\pst@linetype{2}%
\end@ClosedObj}
\def\tx@CCA{CCA }
\def\tx@CCA{CCA }
\def\tx@CC{CC }
\def\tx@IC{IC }
\def\tx@BOC{BOC }
\def\tx@NC{NC }
\def\tx@EOC{EOC }
\def\tx@BAC{BAC }
\def\tx@NAC{NAC }
\def\tx@EAC{EAC }
\def\tx@OpenCurve{OpenCurve }
\def\tx@AltCurve{AltCurve }
\def\tx@ClosedCurve{ClosedCurve }
\def\psset@curvature#1{%
\edef\pst@tempg{#1 }%
\expandafter\psset@@curvature\pst@tempg * * * \@nil}
\def\psset@@curvature#1 #2 #3 #4\@nil{%
\pst@checknum{#1}\pst@tempg
\pst@checknum{#2}\pst@temph
\pst@checknum{#3}\pst@tempi
\edef\psk@curvature{\pst@tempg \pst@temph \pst@tempi}}
\psset@curvature{1 .1 0}
\def\pscurve{\pst@object{pscurve}}
\def\pscurve@i{%
\pst@getarrows{%
\begin@OpenObj
\pst@getcoors[\pscurve@ii}}
\def\pscurve@ii{%
\addto@pscode{%
\pst@cp
\psk@curvature\space /c ED /b ED /a ED
\ifshowpoints true \else false \fi
\tx@OpenCurve}%
\end@OpenObj}
\def\psecurve{\pst@object{psecurve}}
\def\psecurve@i{%
\pst@getarrows{%
\begin@OpenObj
\pst@getcoors[\psecurve@ii}}
\def\psecurve@ii{%
\addto@pscode{%
\psk@curvature\space /c ED /b ED /a ED
\ifshowpoints true \else false \fi
\tx@AltCurve}%
\end@OpenObj}
\def\psccurve{\pst@object{psccurve}}
\def\psccurve@i{%
\begin@ClosedObj
\pst@getcoors[\psccurve@ii}
\def\psccurve@ii{%
\addto@pscode{%
\psk@curvature\space /c ED /b ED /a ED
\ifshowpoints true \else false \fi
\tx@ClosedCurve}%
\def\pst@linetype{1}%
\end@ClosedObj}
\def\psset@dotsize#1{%
\pst@expandafter\pst@getdimnum{#1} 0 {} {}\@nil
\edef\psk@@dotsize{\pst@number\pst@dimg}%
\let\psk@@@dotsize\pst@tempg
\edef\psk@dotsize{%
/DS \psk@@dotsize \psk@@@dotsize CLW mul add 2 div def }}
\psset@dotsize{2pt 2}
\def\psset@dotscale#1{%
\pst@getscale{#1}\psk@dotscale
\ifx\psk@dotscale\@empty
\def\psk@xdotscale{1 }%
\def\psk@ydotscale{1 }%
\else
\let\psk@xdotscale\pst@tempg
\let\psk@ydotscale\pst@temph
\fi}
\def\pst@Getangle#1#2{%
\pst@getangle{#1}\pst@tempg
\def\pst@temph{0. }%
\ifx\pst@tempg\pst@temph
\def#2{}%
\else
\edef#2{\pst@tempg\space rotate }%
\fi}
\def\psset@dotangle#1{%
\pst@getangle{#1}\psk@@dotangle
\ifdim\psk@@dotangle\p@=\z@
\let\psk@dotangle\@empty
\else
\edef\psk@dotangle{\psk@@dotangle rotate }%
\fi}
\psset@dotangle{0}
\def\pst@getdotsize{%
\pst@dimg=\psk@@@dotsize\pslinewidth
\advance\pst@dimg\psk@@dotsize\p@
\pst@dimh=\psk@ydotscale\pst@dimg
\pst@dimg=\psk@xdotscale\pst@dimg
\divide\pst@dimh 2
\divide\pst@dimg 2\relax}
\psset@dotscale{1}
\def\psdot{\pst@object{psdot}}
\def\psdot@i{\@ifnextchar({\psdot@ii}{\psdot@ii(\z@,\z@)}}
\def\psdot@ii(#1){%
\begin@SpecialObj
\pst@@getcoor{#1}%
\addto@pscode{%
\psk@dotsize
\@nameuse{psds@\psk@dotstyle}%
\pst@coor Dot}%
\end@SpecialObj}
\def\psdots{\pst@object{psdots}}
\def\psdots@i{%
\begin@SpecialObj
\pst@getcoors[\psdots@ii}
\def\psdots@ii{%
\addto@pscode{false \tx@NArray \psdots@iii}%
\end@SpecialObj}
\def\psdots@iii{%
\psk@dotsize
\@nameuse{psds@\psk@dotstyle}
newpath
n { transform floor .5 add exch floor .5 add exch itransform Dot } repeat}
\def\tx@SQ{SQ }
\def\tx@ST{ST }
\def\tx@SP{SP }
\def\pst@gdot#1{/Dot { gsave T \psk@dotangle \psk@dotscale #1 grestore } def }
\@namedef{psds@*}{\pst@gdot{0 0 DS \tx@SD}}
\@namedef{psds@o}{%
/r2 DS CLW sub def
\pst@gdot{0 0 DS \tx@SD \pst@usecolor\psfillcolor 0 0 r2 \tx@SD}}
\@namedef{psds@square*}{%
/r1 DS .886 mul def
\pst@gdot{r1 \tx@SQ}}
\@namedef{psds@square}{%
/r1 DS .886 mul def /r2 r1 CLW sub def
\pst@gdot{r1 \tx@SQ \pst@usecolor\psfillcolor r2 \tx@SQ}}
\@namedef{psds@triangle*}{%
/y1 DS .778 mul neg def /x1 y1 1.732 mul neg def
\pst@gdot{x1 y1 \tx@ST}}
\@namedef{psds@triangle}{%
/y1 DS .778 mul neg def /x1 y1 1.732 mul neg def
/y2 y1 CLW add def /x2 y2 1.732 mul neg def
\pst@gdot{x1 y1 \tx@ST \pst@usecolor\psfillcolor x2 y2 \tx@ST}}
\@namedef{psds@pentagon*}{%
/r1 DS 1.149 mul def
\pst@gdot{r1 \tx@SP}}
\@namedef{psds@pentagon}{%
DS .93 mul dup 1.236 mul /r1 ED CLW sub 1.236 mul /r2 ED
\pst@gdot{r1 \tx@SP \pst@usecolor\psfillcolor r2 \tx@SP}}
\@namedef{psds@+}{%
/DS DS 1.253 mul def
\pst@gdot{DS 0 moveto DS neg 0 L stroke 0 DS moveto 0 DS neg L stroke}}
\@namedef{psds@|}{%
\psk@tbarsize CLW mul add 2 div /DS ED
\pst@gdot{0 DS moveto 0 DS neg L stroke}}
\def\psset@dotstyle#1{%
\@ifundefined{psds@#1}%
{\@pstrickserr{Dot style `#1' not defined}\@eha}%
{\edef\psk@dotstyle{#1}}}
\psset@dotstyle{*}
\def\tx@FontDot{FontDot }
\def\newpsfontdot#1[#2]#3#4{%
\@namedef{psds@#1}{%
/#3 \psk@@dotangle [#2] \tx@FontDot
/Dot { moveto gsave \psk@dotscale #4 show grestore } bind def }}
\def\newpsfontdotH#1[#2]#3#4#5{%
\@namedef{psds@#1}{%
/#3 \psk@@dotangle [#2] \tx@FontDot
/Dot {
moveto
\iftrue
gsave \psk@dotscale \pst@usecolor\psfillcolor #5 show grestore
\fi
gsave \psk@dotscale #4 show grestore
} bind def }}
\pstheader{pst-dots.pro}
\newpsfontdot{*}%
[1.0 0.0 0.0 1.0 0.0 0.0]{PSTricksDotFont}{(b)}
\newpsfontdotH{o}%
[1.0 0.0 0.0 1.0 0.0 0.0]{PSTricksDotFont}{(c)}{(b)}
\newpsfontdotH{Bo}%
[1.0 0.0 0.0 1.0 0.0 0.0]{PSTricksDotFont}{(C)}{(b)}
\newpsfontdotH{triangle}%
[1.0 0.0 0.0 1.0 0.0 0.0]{PSTricksDotFont}{(t)}{(u)}
\newpsfontdotH{Btriangle}%
[1.0 0.0 0.0 1.0 0.0 0.0]{PSTricksDotFont}{(T)}{(u)}
\newpsfontdot{triangle*}%
[1.0 0.0 0.0 1.0 0.0 0.0]{PSTricksDotFont}{(u)}
\newpsfontdotH{square}%
[1.0 0.0 0.0 1.0 0.0 0.0]{PSTricksDotFont}{(s)}{(r)}
\newpsfontdotH{Bsquare}%
[1.0 0.0 0.0 1.0 0.0 0.0]{PSTricksDotFont}{(S)}{(r)}
\newpsfontdot{square*}%
[1.0 0.0 0.0 1.0 0.0 0.0]{PSTricksDotFont}{(r)}
\newpsfontdotH{pentagon}%
[1.0 0.0 0.0 1.0 0.0 0.0]{PSTricksDotFont}{(p)}{(q)}
\newpsfontdotH{Bpentagon}%
[1.0 0.0 0.0 1.0 0.0 0.0]{PSTricksDotFont}{(P)}{(q)}
\newpsfontdot{pentagon*}%
[1.0 0.0 0.0 1.0 0.0 0.0]{PSTricksDotFont}{(q)}
\newpsfontdotH{diamond}%
[1.0 0.0 0.0 1.0 0.0 0.0]{PSTricksDotFont}{(d)}{(l)}
\newpsfontdotH{Bdiamond}%
[1.0 0.0 0.0 1.0 0.0 0.0]{PSTricksDotFont}{(D)}{(l)}
\newpsfontdot{diamond*}%
[1.0 0.0 0.0 1.0 0.0 0.0]{PSTricksDotFont}{(l)}
\newpsfontdot{oplus}%
[1.44928 0.0 0.0 1.44928 -0.562319 -0.478261]{Symbol}{<C5>}
\newpsfontdot{otimes}%
[1.44928 0.0 0.0 1.44928 -0.562319 -0.475362]{Symbol}{<C4>}
\newpsfontdot{x}%
[1.8 0.0 0.0 1.8 -0.495 -0.4788]{Symbol}{<B4>}
\newpsfontdot{+}%
[2.3 0.0 0.0 2.3 -0.6486 -0.5819]{Times-Roman}{<2B>}
\newpsfontdot{asterisk}%
[2.43309 0.0 0.0 2.43309 -0.609489 -1.14477]{Times-Roman}{<2A>}
\newpsfontdot{B+}%
[2.3 0.0 0.0 2.3 -0.6555 -0.5819]{Times-Bold}{<2B>}
\newpsfontdot{Basterisk}%
[2.29358 0.0 0.0 2.29358 -0.576835 -1.08486]{Times-Bold}{<2A>}
\newpsfontdot{|}%
[1.98413 0.0 0.0 1.38 -0.258929 -0.5]{Helvetica}{(|)}
\newpsfontdot{B|}%
[1.98413 0.0 0.0 1.38 -0.277778 -0.5]{Helvetica-Bold}{(|)}
\newdimen\pslinearc
\def\psset@linearc#1{\pssetlength\pslinearc{#1}}
\psset@linearc{0pt}
\def\psline{\pst@object{psline}}
\def\psline@i{%
\pst@getarrows{%
\begin@OpenObj
\pst@getcoors[\psline@ii}}
\def\psline@ii{%
\addto@pscode{\pst@cp \psline@iii \tx@Line}%
\end@OpenObj}
\def\psline@iii{%
\ifdim\pslinearc>\z@
/r \pst@number\pslinearc def
/Lineto { \tx@Arcto } def
\else
/Lineto /lineto load def
\fi
\ifshowpoints true \else false \fi}
\def\qline(#1)(#2){%
\def\pst@par{}%
\begin@SpecialObj
\def\pst@linetype{0}%
\pst@getcoor{#1}\pst@tempa
\pst@@getcoor{#2}%
\addto@pscode{%
\pst@tempa moveto \pst@coor L
\@nameuse{psls@\pslinestyle}}%
\end@SpecialObj}
\def\pspolygon{\pst@object{pspolygon}}
\def\pspolygon@i{%
\begin@ClosedObj
\def\pst@cp{}%
\pst@getcoors[\pspolygon@ii}
\def\pspolygon@ii{%
\addto@pscode{\psline@iii \tx@Polygon}%
\def\pst@linetype{1}%
\end@ClosedObj}
\def\psset@framearc#1{\pst@checknum{#1}\psk@framearc}
\psset@framearc{0}
\def\psset@cornersize#1{%
\pst@expandafter\psset@@cornersize{#1}\@nil}
\def\psset@@cornersize#1#2\@nil{%
\if #1a\relax
\def\psk@cornersize{\pst@number\pslinearc false }%
\else
\def\psk@cornersize{\psk@framearc true }%
\fi}
\psset@cornersize{relative}
\def\tx@Rect{Rect }
\def\tx@OvalFrame{OvalFrame }
\def\tx@Frame{Frame }
\def\psset@dimen#1{%
\pst@expandafter\psset@@dimen{#1}\@nil}
\def\psset@@dimen#1#2\@nil{%
\if #1o\relax
\def\psk@dimen{.5 }%
\else
\if #1m\relax
\def\psk@dimen{0 }%
\else
\if #1i\relax
\def\psk@dimen{-.5 }%
\fi
\fi
\fi}
\psset@dimen{outer}
\def\psframe{\pst@object{psframe}}
\def\psframe@i(#1){%
\@ifnextchar({\psframe@ii(#1)}{\psframe@ii(0,0)(#1)}}
\def\psframe@ii(#1)(#2){%
\begin@ClosedObj
\pst@getcoor{#1}\pst@tempa
\pst@@getcoor{#2}%
\addto@pscode{\psk@cornersize \pst@tempa \pst@coor \psk@dimen \tx@Frame}%
\def\pst@linetype{2}%
\showpointsfalse
\end@ClosedObj}
\def\tx@BezierNArray{BezierNArray }
\def\tx@OpenBezier{OpenBezier }
\def\tx@ClosedBezier{ClosedBezier }
\def\tx@BezierShowPoints{BezierShowPoints }
\def\psbezier{\pst@object{psbezier}}
\def\psbezier@i{%
\pst@getarrows{%
\begin@OpenObj
\pst@getcoors[\psbezier@ii}}
\def\psbezier@ii{%
\addto@pscode{%
\ifshowpoints true \else false \fi
\tx@OpenBezier
\ifshowpoints \tx@BezierShowPoints \fi}%
\end@OpenObj}
\def\pscbezier{\pst@object{pscbezier}}
\def\pscbezier@i{%
\begin@ClosedObj
\pst@getcoors[\pscbezier@ii}
\def\pscbezier@ii{%
\addto@pscode{%
\ifshowpoints true \else false \fi
\tx@ClosedBezier
\ifshowpoints \tx@BezierShowPoints \fi}%
\chardef\pst@linetype=1
\end@ClosedObj}
\def\tx@Parab{Parab }
\def\parabola{\pst@object{parabola}}
\def\parabola@i{\pst@getarrows\parabola@ii}
\def\parabola@ii#1(#2)#3(#4){%
\begin@OpenObj
\pst@getcoor{#2}\pst@tempa
\pst@@getcoor{#4}%
\addto@pscode{\pst@tempa \pst@coor \tx@Parab}%
\end@OpenObj}
\def\psset@gridwidth#1{\pst@getlength{#1}\psk@gridwidth}
\psset@gridwidth{.8pt}
\def\psset@griddots#1{%
\pst@cntg=#1\relax
\edef\psk@griddots{\the\pst@cntg}}
\psset@griddots{0}
\def\psset@gridcolor#1{\pst@getcolor{#1}\psgridcolor}
\psset@gridcolor{black}
\def\psset@subgridwidth#1{\pst@getlength{#1}\psk@subgridwidth}
\psset@subgridwidth{.4pt}
\def\psset@subgridcolor#1{\pst@getcolor{#1}\pssubgridcolor}
\psset@subgridcolor{gray}
\def\psset@subgriddots#1{%
\pst@cntg=#1\relax\edef\psk@subgriddots{\the\pst@cntg}}
\psset@subgriddots{0}
\def\psset@subgriddiv#1{%
\pst@cntg=#1\relax\edef\psk@subgriddiv{\the\pst@cntg}}
\psset@subgriddiv{5}
\def\psset@gridlabels#1{\pst@getlength{#1}\psk@gridlabels}
\psset@gridlabels{10pt}
\def\psset@gridlabelcolor#1{\pst@getcolor{#1}\psgridlabelcolor}
\psset@gridlabelcolor{black}
\def\tx@Grid{Grid }
\def\psgrid{\pst@object{psgrid}}
\def\psgrid@i{\@ifnextchar(%
{\psgrid@ii}{\expandafter\psgrid@iv\pic@coor}}
\def\psgrid@ii(#1){\@ifnextchar(%
{\psgrid@iii(#1)}{\psgrid@iv(0,0)(0,0)(#1)}}
\def\psgrid@iii(#1)(#2){\@ifnextchar(%
{\psgrid@iv(#1)(#2)}{\psgrid@iv(#1)(#1)(#2)}}
\def\psgrid@iv(#1)(#2)(#3){%
\begin@SpecialObj
\pst@getcoor{#1}\pst@tempa
\pst@getcoor{#2}\pst@tempb
\pst@@getcoor{#3}%
\ifnum\psk@subgriddiv>1
\addto@pscode{gsave
\psk@subgridwidth SLW \pst@usecolor\pssubgridcolor
\pst@tempb \pst@coor \pst@tempa
\pst@number\psxunit \pst@number\psyunit
\psk@subgriddiv\space \psk@subgriddots\space
{} 0 \tx@Grid grestore}%
\fi
\addto@pscode{gsave
\psk@gridwidth SLW \pst@usecolor\psgridcolor
\pst@tempb \pst@coor \pst@tempa
\pst@number\psxunit \pst@number\psyunit
1 \psk@griddots\space { \pst@usecolor\psgridlabelcolor }
\psk@gridlabels \tx@Grid grestore}%
\end@SpecialObj}
\newif\ifpsmathbox
\psmathboxtrue
\def\pst@mathflag{\z@}
\newtoks\everypsbox
\let\pst@thisbox\relax
\long\def\pst@makenotverbbox#1#2{%
\edef\pst@mathflag{%
\ifpsmathbox\ifmmode\ifinner 1\else 2\fi\else \z@\fi\else \z@\fi}%
\setbox\pst@hbox=\hbox{%
\ifcase\pst@mathflag\or$\m@th\textstyle\or$\m@th\displaystyle\fi
{\pst@thisbox\the\everypsbox#2}%
\ifnum\pst@mathflag>\z@$\fi}%
#1}
\def\pst@makeverbbox#1{%
\def\pst@afterbox{#1}%
\edef\pst@mathflag{%
\ifpsmathbox\ifmmode\ifinner 1\else 2\fi\else \z@\fi\else \z@\fi}%
\afterassignment\pst@beginbox
\setbox\pst@hbox\hbox}
\def\pst@beginbox{%
\ifcase\pst@mathflag\or$\m@th\or$\m@th\displaystyle\fi
\bgroup\aftergroup\pst@endbox
\pst@thisbox
\the\everypsbox}
\def\pst@endbox{%
\ifnum\pst@mathflag>\z@$\fi
\egroup
\pst@afterbox}
\def\pst@makebox{\pst@@makebox}
\def\psverbboxtrue{\def\pst@@makebox{\pst@makeverbbox}}
\def\psverbboxfalse{\def\pst@@makebox{\pst@makenotverbbox}}
\psverbboxfalse
\def\pst@longbox{%
\def\pst@makebox{%
\gdef\pst@makebox{\pst@@makebox}%
\pst@makelongbox}}
\def\pst@makelongbox#1{%
\def\pst@afterbox{#1}%
\edef\pst@mathflag{%
\ifpsmathbox\ifmmode\ifinner 1\else 2\fi\else \z@\fi\else \z@\fi}%
\setbox\pst@hbox\hbox\bgroup
\aftergroup\pst@afterbox
\ifcase\pst@mathflag\or$\m@th\or$\m@th\displaystyle\fi
\begingroup
\pst@thisbox
\the\everypsbox}
\def\pst@endlongbox{%
\endgroup
\ifnum\pst@mathflag>\z@$\fi
\egroup}
\def\pslongbox#1#2{%
\@namedef{#1}{\pst@longbox#2}%
\@namedef{end#1}{\pst@endlongbox}}
\newdimen\psframesep
\def\psset@framesep#1{\pssetlength\psframesep{#1}}
\psset@framesep{3pt}
\newif\ifpsboxsep
\def\psset@boxsep#1{\@nameuse{psboxsep#1}}
\psset@boxsep{true}
\def\pst@useboxpar{%
\use@par
\if@star
\let\pslinecolor\psfillcolor
\solid@star
\let\solid@star\relax
\fi
\ifpsdoubleline \pst@setdoublesep \fi}
\def\psframebox{\pst@object{psframebox}}
\def\psframebox@i{\pst@makebox\psframebox@ii}
\def\psframebox@ii{%
\begingroup
\pst@useboxpar
\pst@dima=\pslinewidth
\advance\pst@dima by \psframesep
\pst@dimc=\wd\pst@hbox\advance\pst@dimc by \pst@dima
\pst@dimb=\dp\pst@hbox\advance\pst@dimb by \pst@dima
\pst@dimd=\ht\pst@hbox\advance\pst@dimd by \pst@dima
\setbox\pst@hbox=\hbox{%
\ifpsboxsep\kern\pst@dima\fi
\begin@ClosedObj
\addto@pscode{%
\psk@cornersize
\pst@number\pst@dima neg
\pst@number\pst@dimb neg
\pst@number\pst@dimc
\pst@number\pst@dimd
.5
\tx@Frame}%
\def\pst@linetype{2}%
\showpointsfalse
\end@ClosedObj
\box\pst@hbox
\ifpsboxsep\kern\pst@dima\fi}%
\ifpsboxsep\dp\pst@hbox=\pst@dimb\ht\pst@hbox=\pst@dimd\fi
\leavevmode\box\pst@hbox
\endgroup}
\def\psdblframebox{\pst@object{psdblframebox}}
\def\psdblframebox@i{\addto@par{doubleline=true}\psframebox@i}
\def\psclip#1{%
\leavevmode
\begingroup
\begin@psclip
\begingroup
\def\use@pscode{%
\pstVerb{%
\pst@dict
/mtrxc CM def
CP CP T
\tx@STV
\psk@origin
\psk@swapaxes
newpath
\pst@code
clip
newpath
mtrxc setmatrix
moveto
0 setgray
end}%
\gdef\pst@code{}}%
\def\@multips(##1)(##2)##3##4{\pst@misplaced\multips}%
\def\nc@object##1##2##3##4{\pst@misplaced{node connection}}%
\hbox to\z@{#1}%
\endgroup
\def\endpsclip{%
\end@psclip
\endgroup}%
\ignorespaces}
\def\endpsclip{\pst@misplaced\endpsclip}
\let\begin@psclip\relax
\def\end@psclip{\pstVerb{currentpoint initclip moveto}}
\def\AltClipMode{%
\def\end@psclip{\pstVerb{\pst@grestore}}%
\def\begin@psclip{\pstVerb{gsave}}}
\def\clipbox{\@ifnextchar[{\clipbox@}{\clipbox@[\z@]}}
\def\clipbox@[#1]{\pst@makebox{\clipbox@@{#1}}}
\def\clipbox@@#1{%
\pssetlength\pst@dimg{#1}%
\leavevmode\hbox{%
\begin@psclip
\pst@Verb{%
CM \tx@STV CP T newpath
/a \pst@number\pst@dimg def
/w \pst@number{\wd\pst@hbox}a add def
/d \pst@number{\dp\pst@hbox}a add neg def
/h \pst@number{\ht\pst@hbox}a add def
a neg d moveto
a neg h L
w h L
w d L
closepath
clip
newpath
0 0 moveto
setmatrix}%
\unhbox\pst@hbox
\end@psclip}}
\def\psshadowbox{\pst@object{psshadowbox}}
\def\psshadowbox@i{\pst@makebox\psshadowbox@ii}
\def\psshadowbox@ii{%
\begingroup
\pst@useboxpar
\psshadowtrue
\psboxseptrue
\def\psk@shadowangle{-45 }%
\setbox\pst@hbox=\hbox{\psframebox@ii}%
\pst@dimh=\psk@shadowsize\p@
\pst@dimh=.7071\pst@dimh
\pst@dimg=\dp\pst@hbox
\advance\pst@dimg\pst@dimh
\dp\pst@hbox=\pst@dimg
\pst@dimg=\wd\pst@hbox
\advance\pst@dimg\pst@dimh
\wd\pst@hbox=\pst@dimg
\leavevmode
\box\pst@hbox
\endgroup}
\def\pscirclebox{\pst@object{pscirclebox}}
\def\pscirclebox@i{\pst@makebox\pscirclebox@ii}
\def\pscirclebox@ii{%
\begingroup
\pst@useboxpar
\setbox\pst@hbox=\hbox{%
\pst@nodehook
\pscirclebox@iii
\box\pst@hbox}%
\ifpsboxsep \pscirclebox@sep \fi
\leavevmode
\box\pst@hbox
\endgroup}
\def\pscirclebox@iii{%
\if@star
\pslinewidth\z@
\pstverb{\pst@dict \tx@STP \pst@usecolor\psfillcolor
newpath \pscirclebox@iv \tx@SD end}%
\else
\begin@ClosedObj
\def\pst@linetype{4}\showpointsfalse
\addto@pscode{%
\pscirclebox@iv CLW 2 div add 0 360 arc closepath}%
\end@ClosedObj
\fi}
\def\pscirclebox@iv{%
\pst@number{\wd\pst@hbox}2 div
\pst@number{\ht\pst@hbox}\pst@number{\dp\pst@hbox}add 2 div
2 copy \pst@number{\dp\pst@hbox}sub 4 2 roll
\tx@Pyth \pst@number\psframesep add }
\def\pscirclebox@sep{%
\pst@dimb=\ht\pst@hbox
\advance\pst@dimb\dp\pst@hbox
\divide\pst@dimb 2
\pst@dima=.5\wd\pst@hbox
\pst@pyth\pst@dima\pst@dimb\pst@dimc
\advance\pst@dimc\pslinewidth
\advance\pst@dimc\psframesep
\advance\pst@dimb-\pst@dimc
\setbox\pst@hbox=\hbox to2\pst@dimc{%
\hss
\vbox{\kern-\pst@dimb\box\pst@hbox}%
\hss}%
\advance\pst@dimb-\dp\pst@hbox
\dp\pst@hbox=-\pst@dimb}
\let\pst@nodehook\relax
\def\psCirclebox{\pst@object{psCirclebox}}
\def\psCirclebox@i{\pst@makebox\psCirclebox@ii}
\def\psCirclebox@ii{%
\begingroup
\pst@useboxpar
\pst@dima=\ht\pst@hbox
\advance\pst@dima\dp\pst@hbox
\divide\pst@dima\tw@
\pssetlength\pst@dimb\psk@radius
\setbox\pst@hbox=\hbox{%
\pst@nodehook
\pscircle(.5\wd\pst@hbox,\pst@dima){\pst@dimb}%
\box\pst@hbox}%
\ifpsboxsep \psCirclebox@sep \fi
\leavevmode
\box\pst@hbox
\endgroup}
\def\psCirclebox@sep{%
\pst@dimc=\pst@dimb
\advance\pst@dimb-\pst@dima
\advance\pst@dima\pst@dimc
\setbox\pst@hbox=\hbox to\tw@\pst@dimc{%
\hss
\vrule width \z@ depth \pst@dimb height \pst@dima
\box\pst@hbox
\hss}}%
\def\psovalbox{\pst@object{psovalbox}}
\def\psovalbox@i{\pst@makebox{\psovalbox@ii}}
\def\psovalbox@ii{%
\begingroup
\pst@useboxpar
\psovalbox@iii
\ifpsboxsep\psovalbox@sep\fi
\leavevmode
\box\pst@hbox
\endgroup}
\def\psovalbox@iii{%
\psovalbox@iv
\setbox\pst@hbox=\hbox{%
\begin@ClosedObj
\addto@pscode{%
0 360
\pst@number\pst@dimc CLW 2 div sub
\pst@number\pst@dimd CLW 2 div sub
\pst@number\pst@dima
\pst@number\pst@dimb
\tx@Ellipse
closepath}%
\def\pst@linetype{2}%
\end@ClosedObj
\unhbox\pst@hbox}}
\def\psovalbox@iv{%
\pst@dimc=\pslinewidth\advance\pst@dimc\psframesep
\pst@dimd=\ht\pst@hbox\advance\pst@dimd\dp\pst@hbox
\pst@dima=.5\wd\pst@hbox
\pst@dimb=.5\pst@dimd\advance\pst@dimb-\dp\pst@hbox
\pst@dimd=.707\pst@dimd
\advance\pst@dimd\pst@dimc
\advance\pst@dimc.707\wd\pst@hbox}
\def\psovalbox@sep{%
\setbox\pst@hbox\hbox to 2\pst@dimc{\hss\unhbox\pst@hbox\hss}%
\pst@dimg=\pst@dimd
\advance\pst@dimg-\pst@dimb
\dp\pst@hbox=\pst@dimg
\advance\pst@dimd\pst@dimb
\ht\pst@hbox=\pst@dimd}
\def\psdiabox{\pst@object{psdiabox}}
\def\psdiabox@i{\pst@makebox{\psdiabox@ii}}
\def\psdiabox@ii{%
\begingroup
\pst@useboxpar
\psdiabox@iii
\ifpsboxsep\psdiabox@sep\fi
\leavevmode
\box\pst@hbox
\endgroup}
\def\psdiabox@iv{%
\pst@dimg=.707\pslinewidth
\advance\pst@dimg.707\psframesep
\pst@dima=\wd\pst@hbox
\divide\pst@dima 2
\pst@dimc=\pst@dima
\advance\pst@dimc\pst@dimg
\pst@dimd=\ht\pst@hbox
\advance\pst@dimd\dp\pst@hbox
\divide\pst@dimd 2
\pst@dimb=\pst@dimd
\advance\pst@dimb-\dp\pst@hbox
\advance\pst@dimd\pst@dimg}
\def\psdiabox@iii{%
\psdiabox@iv
\setbox\pst@hbox=\hbox{%
\begin@ClosedObj
\addto@pscode{%
\psline@iii
pop
.5
\pst@number\pst@dimc 2 mul \pst@number\pst@dimd 2 mul
0
\pst@number\pst@dima \pst@number\pst@dimb
\tx@Diamond}%
\def\pst@linetype{4}%
\end@ClosedObj
\box\pst@hbox}}
\def\psdiabox@sep{%
\setbox\pst@hbox\hbox to 4\pst@dimc{\hss\unhbox\pst@hbox\hss}%
\multiply\pst@dimd 2
\advance\pst@dimd\pst@dimb
\ht\pst@hbox\pst@dimd
\advance\pst@dimd-2\pst@dimb
\dp\pst@hbox\pst@dimd}
\def\psset@trimode#1{\pst@expandafter\psset@@trimode{#1}\@empty\@empty\@nil}
\def\psset@@trimode#1#2#3\@nil{%
\let\pst@tempg#1\relax
\ifx\pst@tempg*%
\let\psk@@trimode\@empty
\let\pst@tempg#2\relax
\else
\let\psk@@trimode\relax
\fi
\edef\psk@trimode{%
\ifx R\pst@tempg 1 \else\ifx D\pst@tempg 2
\else\ifx L\pst@tempg 3 \else 0 \fi\fi\fi}}
\psset@trimode{U}
\def\pstribox{\pst@object{pstribox}}
\def\pstribox@i{\pst@makebox{\pstribox@ii}}
\def\pstribox@ii{%
\begingroup
\pst@useboxpar
\pstribox@iii
\ifpsboxsep\pstribox@sep\fi
\leavevmode
\box\pst@hbox
\endgroup}
\def\pstribox@iii{%
\pstribox@iv
\setbox\pst@hbox=\hbox{%
\begin@ClosedObj
\addto@pscode{%
\psline@iii
pop
.5
\pst@number\pst@dimc \pst@number\pst@dimd
\ifodd\psk@trimode exch \fi
\psk@trimode -90 mul
\pst@number\pst@dima \pst@number\pst@dimb
\tx@Triangle}%
\def\pst@linetype{2}%
\end@ClosedObj
\box\pst@hbox}}
\def\pstribox@iv{%
\pst@dimh=\pslinewidth
\advance\pst@dimh\psframesep
\pst@dimg=\ht\pst@hbox
\advance\pst@dimg-\dp\pst@hbox
\divide\pst@dimg 2
\edef\pst@tempa{\number\pst@dimg sp}
\ifodd\psk@trimode
\pst@dimb\pst@dimg
\else
\pst@dima=\wd\pst@hbox
\divide\pst@dima 2
\fi
\ifcase\psk@trimode
\pst@dimb=-\dp\pst@hbox
\advance\pst@dimb-\pst@dimh
\or
\pst@dima=-\pst@dimh
\or
\pst@dimb=\ht\pst@hbox
\advance\pst@dimb\pst@dimh
\or
\pst@dima=\wd\pst@hbox
\advance\pst@dima\pst@dimh
\fi
\pst@dimd=\dp\pst@hbox
\advance\pst@dimd\ht\pst@hbox
\ifx\psk@@trimode\relax
\pst@dimc=\wd\pst@hbox
\advance\pst@dimc\ifodd\psk@trimode 1.447\else 1.789\fi\pst@dimh
\multiply\pst@dimc 2
\advance\pst@dimd\ifodd\psk@trimode 1.789\else 1.447\fi\pst@dimh
\multiply\pst@dimd 2
\else
\ifodd\psk@trimode
\advance\pst@dimd 1.1547\wd\pst@hbox
\advance\pst@dimd 3.4641\pst@dimh
\pst@dimc=.866\pst@dimd
\else
\advance\pst@dimd .866\wd\pst@hbox 
\advance\pst@dimd 3\pst@dimh
\pst@dimc=1.1547\pst@dimd 
\fi
\fi}
\def\pstribox@sep{%
\ifodd\psk@trimode
\advance\pst@dimb.5\pst@dimd
\ht\pst@hbox=\pst@dimb
\advance\pst@dimd-\pst@dimb
\dp\pst@hbox=\pst@dimd
\else
\setbox\pst@hbox\hbox to \pst@dimc{\hss\unhbox\pst@hbox\hss}%
\global\pst@dimg=.5\pst@dimc
\fi
\ifcase\psk@trimode
\dp\pst@hbox-\pst@dimb
\advance\pst@dimd\pst@dimb
\ht\pst@hbox\pst@dimd
\or
\pst@dimg=.5\wd\pst@hbox
\global\advance\pst@dimg-\pst@dima
\setbox\pst@hbox\hbox to \pst@dimc{\kern-\pst@dima\box\pst@hbox\hss}%
\or
\ht\pst@hbox\pst@dimb
\advance\pst@dimd-\pst@dimb
\dp\pst@hbox\pst@dimd
\or
\pst@dimg=\pst@dimc
\advance\pst@dimg-\pst@dima
\global\advance\pst@dimg.5\wd\pst@hbox
\setbox\pst@hbox\hbox to \pst@dimc{%
\hss\box\pst@hbox\kern\psframesep\kern\pslinewidth}%
\fi}
\def\psset@arcsepA#1{\pst@getlength{#1}\psk@arcsepA}
\def\psset@arcsepB#1{\pst@getlength{#1}\psk@arcsepB}
\def\psset@arcsep#1{%
\psset@arcsepA{#1}\let\psk@arcsepB\psk@arcsepA}
\psset@arcsep{0}
\def\tx@ArcArrow{ArcArrow }
\def\psarc{\pst@object{psarc}}
\def\psarc@i{%
\@ifnextchar({\psarc@iii}{\psarc@ii}}
\def\psarc@ii#1{\addto@par{arrows=#1}%
\@ifnextchar({\psarc@iii}{\psarc@iii(0,0)}}
\def\psarc@iii(#1)#2#3#4{%
\begin@OpenObj
\pst@getangle{#3}\pst@tempa
\pst@getangle{#4}\pst@tempb
\pst@@getcoor{#1}%
\pssetlength\pst@dima{#2}%
\addto@pscode{\psarc@iv \psarc@v}%
\gdef\psarc@type{0}%
\showpointsfalse
\end@OpenObj}
\def\psarc@iv{%
\pst@coor /y ED /x ED
/r \pst@number\pst@dima def
/c 57.2957 r \tx@Div def
/angleA
\pst@tempa
\psk@arcsepA c mul 2 div
\ifcase \psarc@type add \or sub \fi
def
/angleB
\pst@tempb
\psk@arcsepB c mul 2 div
\ifcase \psarc@type sub \or add \fi
def
\ifshowpoints\psarc@showpoints\fi
\ifx\psk@arrowA\@empty
\ifnum\psk@liftpen=2
r angleA \tx@PtoC
y add exch x add exch
moveto
\fi
\fi}
\def\psarc@v{%
x y r
angleA
\ifx\psk@arrowA\@empty\else
{ ArrowA CP }
{ \ifcase\psarc@type add \or sub \fi }
\tx@ArcArrow
\fi
angleB
\ifx\psk@arrowB\@empty\else
{ ArrowB }
{ \ifcase\psarc@type sub \or add \fi }
\tx@ArcArrow
\fi
\ifcase\psarc@type arc \or arcn \fi}
\def\psarc@type{0}
\def\psarc@showpoints{%
gsave
newpath
x y moveto
x y r \pst@tempa \pst@tempb
\ifcase\psarc@type arc \or arcn \fi
closepath
CLW 2 div SLW
[ \psk@dash\space ] 0 setdash stroke
grestore }
\def\psarcn{\pst@object{psarcn}}
\def\psarcn@i{\def\psarc@type{1}\psarc@i}
\def\pscircle{\pst@object{pscircle}}
\def\pscircle@i{\@ifnextchar({\pscircle@do}{\pscircle@do(0,0)}}
\def\pscircle@do(#1)#2{%
\if@star
{\use@par\qdisk(#1){#2}}%
\else
\begin@ClosedObj
\pst@@getcoor{#1}%
\pssetlength\pst@dimc{#2}%
\def\pst@linetype{4}%
\addto@pscode{%
\pst@coor
\pst@number\pst@dimc
\psk@dimen CLW mul sub
0 360 arc
closepath}%
\showpointsfalse
\end@ClosedObj
\fi
\ignorespaces}
\def\qdisk(#1)#2{%
\def\pst@par{}%
\begin@SpecialObj
\pst@@getcoor{#1}%
\pssetlength\pst@dimg{#2}%
\addto@pscode{\pst@coor \pst@number\pst@dimg \tx@SD}%
\end@SpecialObj}
\def\psset@radius#1{\pst@@getlength{#1}\psk@radius}
\psset@radius{.25cm}
\def\psCircle{\pst@object{psCircle}}
\def\psCircle@i{\@ifnextchar({\psCircle@ii}{\psCircle@ii(0,0)}}
\def\psCircle@ii(#1){\pscircle@do(#1){\psk@radius}}
\def\pswedge{\pst@object{pswedge}}
\def\pswedge@i{\@ifnextchar({\pswedge@ii}{\pswedge@ii(0,0)}}
\def\pswedge@ii(#1)#2#3#4{%
\begin@ClosedObj
\pssetlength\pst@dimc{#2}
\pst@getangle{#3}\pst@tempa
\pst@getangle{#4}\pst@tempb
\pst@@getcoor{#1}%
\def\pst@linetype{1}%
\addto@pscode{%
\pst@coor
2 copy
moveto
\pst@number\pst@dimc \psk@dimen CLW mul sub 
\pst@tempa \pst@tempb
arc
closepath}%
\showpointsfalse
\end@ClosedObj}
\def\tx@Ellipse{Ellipse }
\def\psellipse{\pst@object{psellipse}}
\def\psellipse@i(#1){\@ifnextchar(%
{\psellipse@ii(#1)}{\psellipse@ii(0,0)(#1)}}
\def\psellipse@ii(#1)(#2){%
\begin@ClosedObj
\pst@getcoor{#1}\pst@tempa
\pst@@getcoor{#2}%
\addto@pscode{%
0 360
\pst@coor
\ifdim\psk@dimen\p@=\z@\else
\psk@dimen CLW mul
dup 4 -1 roll sub neg 3 1 roll sub
\fi
\pst@tempa
\tx@Ellipse
closepath}%
\def\pst@linetype{2}%
\end@ClosedObj}
\def\multips{\@ifnextchar({\def\pst@par{}\multips@ii}{\multips@i}}
\def\multips@i#1{\def\pst@par{rot=#1}\multips@ii}
\def\multips@ii(#1){\@ifnextchar({\multips@iii(#1)}{\multips@iii(\z@,\z@)(#1)}}
\long\def\multips@iii(#1)(#2)#3#4{%
\begingroup
\use@par
\pst@getcoor{#1}\pst@tempa
\pst@@getcoor{#2}%
\pst@cnta=#3\relax
\init@pscode
\addto@pscode{%
\pst@tempa T \the\pst@cnta\space \pslbrace
gsave \ifx\psk@rot\@empty\else\psk@rot rotate \fi}%
\hbox to\z@{%
\def\init@pscode{%
\addto@pscode{%
gsave
\pst@number\pslinewidth SLW
\pst@usecolor\pslinecolor}}%
\def\use@pscode{\addto@pscode{grestore}}%
\def\psclip##1{\pst@misplaced\psclip}%
\def\nc@object##1##2##3##4{\pst@misplaced{node connection}}%
#4}%
\addto@pscode{grestore \pst@coor T \psrbrace repeat}%
\leavevmode
\use@pscode
\endgroup
\ignorespaces}
\def\scalebox#1{\pst@makebox{\@scalebox{#1}}}
\def\@scalebox#1{%
\begingroup
\pst@getscale{#1}\pst@tempa
\let\pst@tempc\pst@tempg
\let\pst@tempd\pst@temph
\@@scalebox
\endgroup}
\def\@@scalebox{%
\leavevmode
\hbox{%
\ifdim\pst@tempd\p@<\z@
\pst@dimg=\pst@tempd\ht\pst@hbox
\pst@dimh=\pst@tempd\dp\pst@hbox
\dp\pst@hbox=-\pst@dimg
\ht\pst@hbox=-\pst@dimh
\else
\ht\pst@hbox=\pst@tempd\ht\pst@hbox
\dp\pst@hbox=\pst@tempd\dp\pst@hbox
\fi
\pst@dima=\pst@tempc\wd\pst@hbox
\ifdim\pst@dima<\z@\kern-\pst@dima\fi
\pst@Verb{CP CP translate \pst@tempa \tx@NET}%
\hbox to \z@{\box\pst@hbox\hss}%
\pst@Verb{%
CP CP translate
1 \pst@tempc div 1 \pst@tempd div scale
\tx@NET}%
\ifdim\pst@dima>\z@\kern\pst@dima\fi}}
\pslongbox{Scalebox}{\scalebox}
\def\scaleboxto(#1,#2){\pst@makebox{\@scaleboxto(#1,#2)}}
\def\@scaleboxto(#1,#2){%
\begingroup
\pssetlength\pst@dima{#1}%
\pssetlength\pst@dimb{#2}%
\ifdim\pst@dima=\z@\else
\pst@divide{\pst@dima}{\wd\pst@hbox}\pst@tempc
\edef\pst@tempc{\pst@tempc\space}%
\fi
\ifdim\pst@dimb=\z@
\ifdim\pst@dima=\z@
\@pstrickserr{%
\string\scaleboxto\space dimensions cannot both be zero}\@ehpa
\def\pst@tempa{}%
\def\pst@tempc{1 }%
\def\pst@tempd{1 }%
\else
\let\pst@tempd\pst@tempc
\fi
\else
\pst@dimc=\ht\pst@hbox
\advance\pst@dimc\dp\pst@hbox
\pst@divide{\pst@dimb}{\pst@dimc}\pst@tempd
\edef\pst@tempd{\pst@tempd\space}%
\ifdim\pst@dima=\z@ \let\pst@tempc\pst@tempd \fi
\fi
\edef\pst@tempa{\pst@tempc \pst@tempd scale }%
\@@scalebox
\endgroup}
\pslongbox{Scaleboxto}{\scaleboxto}
\def\tx@Rot{Rot }
\def\rotateleft{\pst@makebox{\@rotateleft\pst@hbox}}
\def\@rotateleft#1{%
\leavevmode\hbox{\hskip\ht#1\hskip\dp#1\vbox{\vskip\wd#1%
\pst@Verb{90 \tx@Rot}
\vbox to \z@{\vss\hbox to \z@{\box#1\hss}\vskip\z@}%
\pst@Verb{-90 \tx@Rot}}}}
\def\rotateright{\pst@makebox{\@rotateright\pst@hbox}}
\def\@rotateright#1{%
\hbox{\hskip\ht#1\hskip\dp#1\vbox{\vskip\wd#1%
\pst@Verb{-90 \tx@Rot}
\vbox to \z@{\hbox to \z@{\hss\box#1}\vss}%
\pst@Verb{90 \tx@Rot}}}}
\def\rotatedown{\pst@makebox{\@rotatedown\pst@hbox}}
\def\@rotatedown#1{%
\hbox{\hskip\wd#1\vbox{\vskip\ht#1\vskip\dp#1%
\pst@Verb{180 \tx@Rot}%
\vbox to \z@{\hbox to \z@{\box#1\hss}\vss}%
\pst@Verb{-180 \tx@Rot}}}}
\pslongbox{Rotateleft}{\rotateleft}
\pslongbox{Rotateright}{\rotateright}
\pslongbox{Rotatedown}{\rotatedown}
\def\pst@starbox{%
\setbox\pst@hbox\hbox{\psframebox*[boxsep=false]{\unhbox\pst@hbox}}}
\def\pst@@makesmall#1{%
\setbox#1=\hbox to\z@{\hss\vbox to \z@{\vss\box#1\vss}\hss}}
\def\pst@@@makesmall#1{%
\pst@dimh=\psk@xref\wd#1%
\ifx\psk@yref\relax
\pst@dimg=\dp#1%
\else
\pst@dimg=\psk@yref\ht#1%
\advance\pst@dimg\psk@yref\dp#1%
\fi
\setbox#1=\hbox to\z@{%
\kern-\pst@dimh\vbox to\z@{\vss\box#1\kern-\pst@dimg}\hss}}
\def\psset@ref#1{\pst@expandafter\psset@@ref{#1}\@empty,,\@nil}
\def\psset@@ref#1#2,#3,#4\@nil{%
\def\psk@xref{.5}%
\def\psk@yref{.5}%
\let\pst@makesmall\pst@@@makesmall
\ifx\@empty#3\@empty
\@nameuse{getref@#1}%
\@nameuse{getref@#2}%
\else
\pst@checknum{#1#2}\psk@xref
\pst@checknum{#3}\psk@yref
\fi}
\def\getref@c{\let\pst@makesmall\pst@@makesmall}
\def\getref@t{\def\psk@yref{1}}
\def\getref@b{\def\psk@yref{0}}
\def\getref@B{\let\psk@yref\relax}
\def\getref@l{\def\psk@xref{0}}
\def\getref@r{\def\psk@xref{1}}
\psset@ref{c}
\def\psset@rot#1{%
\pst@expandafter{\@ifnextchar*{\psset@@@rot}{\psset@@rot}}{#1}\@nil}
\def\psset@@rot#1\@nil{%
\def\next##1@#1=##2@##3\@nil{%
\ifx\relax##2%
\pst@getangle{#1}\psk@rot
\else
\def\psk@rot{##2}%
\fi}%
\expandafter\next\pst@rottable @#1=\relax @\@nil}
\def\psset@@@rot#1#2\@nil{%
\psset@@rot#2\@nil
\edef\psk@rot{\pst@rotlist \ifx\psk@rot\@empty\else\psk@rot add \fi}}
\def\pst@rotlist{mark RAngle /a ED cleartomark a neg }
\def\pst@rottable{%
@0=%
@U=%
@L=90 %
@D=180 %
@R=-90 %
@N=\pst@rotlist
@W=\pst@rotlist 90 add %
@S=\pst@rotlist 180 add %
@E=\pst@rotlist 90 sub }
\psset@rot{0}
\def\tx@RotBegin{RotBegin }
\def\tx@RotEnd{RotEnd }
\def\pst@rotate#1#2{%
\ifx#1\@empty\else
\setbox#2=\hbox{\pst@Verb{#1 \tx@RotBegin}\box#2\pst@Verb{\tx@RotEnd}}%
\fi}
\def\psput@cartesian#1{%
\hbox to \z@{\kern\pst@dimg{\vbox to \z@{\vss\box#1\vskip\pst@dimh}\hss}}}
\def\psput@special#1{%
\hbox{%
\pst@Verb{{ \pst@coor } \tx@PutCoor \tx@PutBegin}%
\box#1%
\pst@Verb{\tx@PutEnd}}}
\def\tx@PutCoor{PutCoor }
\def\tx@PutBegin{PutBegin }
\def\tx@PutEnd{PutEnd }
\def\rput{\def\pst@par{}\pst@ifstar{\@ifnextchar[{\rput@i}{\rput@ii}}}
\def\rput@i[#1]{\addto@par{ref={#1}}\rput@ii}
\def\rput@ii{\@ifnextchar({\rput@iv}{\rput@iii}}
\def\rput@iii#1{\addto@par{rot={#1}}\@ifnextchar({\rput@iv}{\rput@iv(\z@,\z@)}}
\def\rput@iv(#1){\pst@killglue\pst@makebox{\rput@v{#1}}}
\def\rput@v#1{%
\begingroup
\use@par
\if@star\pst@starbox\fi
\pst@makesmall\pst@hbox
\pst@rotate\psk@rot\pst@hbox
\psput@{#1}\pst@hbox
\endgroup
\ignorespaces}
\def\multirput{%
\def\pst@par{}%
\pst@ifstar{\@ifnextchar[{\multirput@i}{\multirput@ii}}}
\def\multirput@i[#1]{\addto@par{ref={#1}}\multirput@ii}
\def\multirput@ii{\@ifnextchar({\multirput@iv}{\multirput@iii}}
\def\multirput@iii#1{\addto@par{rot={#1}}\multirput@iv}
\def\multirput@iv(#1){%
\@ifnextchar({\multirput@v(#1)}{\multirput@v(\z@,\z@)(#1)}}
\def\multirput@v(#1,#2)(#3,#4)#5{%
\pst@makebox{\multirput@vi(#1,#2)(#3,#4){#5}}}
\def\multirput@vi(#1,#2)(#3,#4)#5{%
\begingroup
\use@par
\if@star\pst@starbox\fi
\pst@makesmall\pst@hbox
\pst@rotate\psk@rot\pst@hbox
\pssetxlength\pst@dima{#1}%
\pssetylength\pst@dimb{#2}%
\pssetxlength\pst@dimc{#3}%
\pssetylength\pst@dimd{#4}%
\pst@cntg=#5\relax
\pst@cnth=\@ne
\leavevmode
\loop
\vbox to \z@{%
\vss
\hbox to \z@{\kern\pst@dima\copy\pst@hbox\hss}%
\vskip\pst@dimb}%
\ifnum\pst@cntg>\pst@cnth
\advance\pst@dima\pst@dimc
\advance\pst@dimb\pst@dimd
\advance\pst@cnth\@ne
\repeat
\endgroup
\ignorespaces}
\newif\if@fixedradius
\def\cput{\pst@object{cput}}
\def\cput@i{\@fixedradiusfalse\cput@ii}
\def\cput@ii{\pst@killglue\@ifnextchar({\cput@iv}{\cput@iii}}
\def\cput@iii#1{%
\addto@par{rot={#1}}%
\@ifnextchar({\cput@iv}{\cput@iv(\z@,\z@)}}
\def\cput@iv(#1){\pst@makebox{\cput@v{#1}}}
\def\cput@v#1{%
\begingroup
\use@par
\setbox\pst@hbox=\hbox{%
\psboxsepfalse
\if@fixedradius\psCirclebox@ii\else\pscirclebox@ii\fi}%
\pst@@makesmall\pst@hbox
\pst@rotate\psk@rot\pst@hbox
\psput@{#1}\pst@hbox
\endgroup
\ignorespaces}
\def\Cput{\pst@object{Cput}}
\def\Cput@i{\@fixedradiustrue\cput@ii}
\newdimen\pslabelsep
\def\psset@labelsep#1{\pssetlength\pslabelsep{#1}}
\psset@labelsep{5pt}
\def\psset@refangle#1{\pst@expandafter\psset@@refangle{#1}\@nil}
\def\psset@@refangle#1\@nil{%
\def\next##1@#1=##2"##3@##4\@nil{%
\ifx\relax##2%
\pst@getangle{#1}\psk@refangle
\def\psk@uputref{}%
\else
\def\psk@refangle{##2 }%
\def\psk@uputref{##3}%
\fi}%
\expandafter\next\pst@refangletable @#1=\relax"@\@nil}
\def\pst@refangletable{%
@r=0"20%
@u=90"02%
@l=180"10%
@d=-90"01%
@ur=45"22%
@ul=135"12%
@dr=-135"21%
@dl=-45"11}
\psset@refangle{0}
\def\uput{\def\pst@par{}\pst@ifstar{\@ifnextchar[{\uput@ii}{\uput@i}}} 
\def\uput@i#1{\addto@par{labelsep=#1}\uput@ii}
\def\uput@ii[#1]{%
\addto@par{refangle={#1}}%
\@ifnextchar({\uput@iv}{\uput@iii}}
\def\uput@iii#1{%
\addto@par{rot={#1}}%
\@ifnextchar({\uput@iv}{\uput@iv(\z@,\z@)}}
\def\uput@iv(#1){\pst@killglue\pst@makebox{\uput@v{#1}}}
\def\uput@v#1{%
\begingroup
\use@par
\if@star\pst@starbox\fi
\uput@vi
\psput@{#1}\pst@hbox
\endgroup
\ignorespaces}
\def\uput@vi{%
\ifx\psk@uputref\@empty
\uput@vii\tx@UUput{}%
\else
\ifx\psk@rot\@empty
\expandafter\uput@viii\psk@uputref
\else
\uput@vii\tx@UUput{}%
\fi
\fi}
\def\uput@vii#1#2{%
\edef\pst@coor{%
\pst@number\pslabelsep
#2%
\pst@number{\wd\pst@hbox}%
\pst@number{\ht\pst@hbox}%
\pst@number{\dp\pst@hbox}%
\psk@refangle\space \ifx\psk@rot\@empty\else\psk@rot\space sub \fi
\tx@Uput #1}%
\setbox\pst@hbox=\hbox to\z@{\hss\vbox to\z@{\vss\box\pst@hbox\vss}\hss}%
\setbox\pst@hbox=\psput@special\pst@hbox
\ifx\psk@rot\@empty\else\pst@rotate\psk@rot\pst@hbox\fi}
\def\uput@viii#1#2{%
\ifnum#1>\z@\ifnum#2>\z@\pslabelsep=.707\pslabelsep\fi\fi
\setbox\pst@hbox=\vbox to\z@{%
\ifnum#2=1 \vskip\pslabelsep\else\vss\fi
\hbox to\z@{%
\ifnum#1=2 \hskip\pslabelsep\else\hss\fi
\box\pst@hbox
\ifnum#1=1 \hskip\pslabelsep\else\hss\fi}%
\ifnum#2=2 \vskip\pslabelsep\else\vss\fi}}
\def\tx@Uput{Uput }
\def\tx@UUput{UUput }
\def\Rput{\def\pst@par{}\pst@ifstar{\@ifnextchar[{\Rput@ii}{\Rput@i}}}
\def\Rput@i#1{\addto@par{labelsep=#1}\Rput@ii}
\def\Rput@ii[#1]{\addto@par{ref={#1}}\@ifnextchar({\Rput@iv}{\Rput@iii}}
\def\Rput@iii#1{\addto@par{rot={#1}}\@ifnextchar({\Rput@iv}{\Rput@iv(\z@,\z@)}}
\def\Rput@iv(#1){\pst@killglue\pst@makebox{\Rput@v{#1}}}
\def\Rput@v#1{%
\begingroup
\use@par
\if@star\pst@starbox\fi
\Rput@vi
\pst@makesmall\pst@hbox
\pst@rotate\psk@rot\pst@hbox
\psput@{#1}\pst@hbox
\endgroup
\ignorespaces}
\def\Rput@vi{%
\pst@dimg=\dp\pst@hbox
\advance\pst@dimg\pslabelsep
\dp\pst@hbox=\pst@dimg
\pst@dimg=\ht\pst@hbox
\advance\pst@dimg\pslabelsep
\ht\pst@hbox=\pst@dimg
\setbox\pst@hbox\hbox{\kern\pslabelsep\box\pst@hbox\kern\pslabelsep}}%
\def\oldpsput{%
\def\pst@par{}\pst@ifstar{\@ifnextchar[{\oldpsput@i}{\oldpsput@ii}}}
\def\oldpsput@i[#1]{\addto@par{ref={#1}}\oldpsput@ii}
\def\oldpsput@ii{\@ifnextchar<{\oldpsput@iii}{\oldpsput@iv}}
\def\oldpsput@iii<#1>{\rput@iii{#1}}
\def\OldPsput{\let\psput\oldpsput}
\def\NewPsput{\let\psput\rput}
\def\pspicture{\begingroup\pst@ifstar\pst@picture}
\def\pst@picture{%
\@ifnextchar[{\pst@@picture}{\pst@@picture[0]}}
\def\pst@@picture[#1]#2(#3,#4){%
\@ifnextchar({\pst@@@picture[#1](#3,#4)}%
{\pst@@@picture[#1](0,0)(#3,#4)}}
\def\pst@@@picture[#1](#2,#3)(#4,#5){%
\pssetxlength\pst@dima{#2}%
\pssetylength\pst@dimb{#3}%
\pssetxlength\pst@dimc{#4}%
\pssetylength\pst@dimd{#5}%
\ifdim\pst@dima>\pst@dimc
\pst@dimg=\pst@dima
\pst@dima=\pst@dimc
\pst@dimc=\pst@dimg
\fi
\ifdim\pst@dimb>\pst@dimd
\pst@dimg=\pst@dimb
\pst@dimb=\pst@dimd
\pst@dimd=\pst@dimg
\fi
\def\pst@tempa{#1}%
\setbox\pst@hbox=\hbox\bgroup
\begingroup\KillGlue
\@ifundefined{@latexerr}{}{\let\unitlength\psunit}%
\edef\pic@coor{(#2,#3)(#2,#3)(#4,#5)}\ignorespaces}
\def\pic@coor{(0,0)(0,0)(10,10)}
\def\endpspicture{%
\pst@killglue
\endgroup
\egroup
\ifdim\wd\pst@hbox=\z@\else
\fi
\ht\pst@hbox=\pst@dimd
\dp\pst@hbox=-\pst@dimb
\setbox\pst@hbox=\hbox{%
\kern-\pst@dima
\ifx\pst@tempa\@empty\else
\advance\pst@dimd-\pst@dimb
\pst@dimd=\pst@tempa\pst@dimd
\advance\pst@dimd\pst@dimb
\lower\pst@dimd
\fi
\box\pst@hbox
\kern\pst@dimc}%
\if@star\setbox\pst@hbox=\hbox{\clipbox@@\z@}\fi
\leavevmode\box\pst@hbox
\endgroup}
\@namedef{pspicture*}{\pspicture*}
\@namedef{endpspicture*}{\endpspicture}
\def\tx@BeginOL{BeginOL }
\def\tx@InitOL{InitOL }
\def\pst@initoverlay#1{\pst@Verb{\tx@InitOL /TheOL (#1) def}}
\def\AltOverlayMode{%
\def\pst@initoverlay##1{%
\pst@Verb{%
\tx@InitOL
/Visible { initclip } def
/Invisible {
CP newpath OLUnit itransform moveto clip newpath moveto
} def
/TheOL (##1) def}}}
\def\pst@overlay#1{%
\edef\curr@overlay{#1}%
\pst@Verb{(#1) BOL}%
\aftergroup\pst@endoverlay}
\def\pst@endoverlay{%
\pst@Verb{(\curr@overlay) BOL}}
\def\curr@overlay{all}
\newbox\theoverlaybox
\def\overlaybox{%
\global\setbox\theoverlaybox=\hbox\bgroup
\begingroup
\let\psoverlay\pst@overlay
\def\overlaybox{%
\@pstrickserr{Overlays cannot be nested}\@eha}%
\def\putoverlaybox{%
\@pstrickserr{You must end the overlay box
before using \string\putoverlaybox}}%
\psoverlay{main}%
\ignorespaces}
\def\endoverlaybox{\endgroup\egroup}
\def\putoverlaybox#1{%
\hbox{\pst@initoverlay{#1}\copy\theoverlaybox}}
\def\psoverlay{\@pstrickserr{\string\psoverlay\space
can only be used after \string\overlaybox}}
\ifx\pstcustomize\relax \input pstricks.con \fi
\catcode`\@=\PstAtCode\relax
\begin{document}
\maketitle

\newtheorem{theo}{Theorem}[section]
\newtheorem{defi}[theo]{Definition}
\newtheorem{defis}[theo]{Definitions}
\newtheorem{rem}[theo]{Remark}
\newtheorem{rems}[theo]{Remarks}
\newtheorem{nota}[theo]{Notation}
\newtheorem{notas}[theo]{Notations}
\newtheorem{prop}[theo]{Proposition}
\newtheorem{hypo}[theo]{Hypothese}
\newtheorem{lem}[theo]{Lemma}
\newtheorem{conv}[theo]{Convention}
\newtheorem{convs}[theo]{Conventions}
\newtheorem{cor}[theo]{Corollary}
\newtheorem{obs}[theo]{Observation}
\newtheorem{obss}[theo]{Observations}
\newtheorem{rap}[theo]{Recall}
\newtheorem{raps}[theo]{Rappels}
\newtheorem{crit}[theo]{Crit\`ere}

$ $
\def\abstractname{Abstract}
\begin{abstract}
$ $

Let $\mathbb{K}$ be a (commutative) field and consider a nonzero element $q$ in $\mathbb{K}$ which is not a root of unity. 
In \cite{glen1}, Goodearl and Lenagan have shown that the number of $\hc$-primes in $R=\mnk$ 
which contain all $(t+1) \times (t+1)$ quantum minors but not all $t \times t$ quantum minors 
is a perfect square. The aim of this paper is to make precise their result: we prove that this number 
is equal to $(t!)^2 S(n+1,t+1)^2$, where $S(n+1,t+1)$ denotes the Stirling number 
of second kind associated to $n+1$ and $t+1$. This result was conjectured by Goodearl, Lenagan and 
McCammond. The proof involves some closed formulas for the poly-Bernoulli numbers that were established 
in \cite{kaneko1} and \cite{kaneko2}.
\end{abstract}
$ $
\\$ $

\section{Introduction.}
$ $

Fix a (commutative) field $\mathbb{K}$ and an integer $n$ greater than or equal to $2$, and choose
 an element $q$ in $\mathbb{K}^*:=\mathbb{K} \setminus\{0\}$ which is not a root of unity. Denote 
by $R=\mnk$ the quantization of the ring of regular functions on 
$n \times n$ matrices with entries in $\mathbb{K}$ and by $(Y_{\ia})_{(\ia) \in \gc 1,n \dc^2}$ 
the matrix of its canonical generators. 
The bialgebra structure of $R$ gives us an action of the group 
$\hc:=(\mathbb{C}^*)^{2n}$ on $R$ by $\mathbb{K}$-automorphisms (See \cite{glen1}) via:
$$(a_1,\dots,a_n,b_1,\dots,b_n).Y_{\ia}=a_i b_{\alpha} Y_{\ia} 
\quad ((\ia) \in \gc 1,n \dc^2).$$

In \cite{gl1}, Goodearl and Letzter have shown that $R$ has only finitely many 
$\hc$-invariant prime ideals (See \cite{gl1}, 5.7. (i)) and that, in order to calculate the prime and primitive spectra of
 $R$, it is enough to determine the $\hc$-invariant prime ideals of $R$ (See \cite{gl1}, Theorem 6.6). 
Next, using the theory of deleting derivations, Cauchon has found a formula for the exact number of 
$\hc$-invariant prime ideals in R (See \cite{c2}, Proprosition 3.3.2). In this paper, we investigate these ideals.

In \cite{lau} (See also \cite{lau3}), we have proved, assuming that $\mathbb{K}=\mathbb{C}$ (the field of complex numbers)
 and $q$ is transcendental over $\mathbb{Q}$, that the $\hc$-invariant prime ideals in $\mnc$ 
are generated by quantum minors, 
as conjectured by Goodearl and Lenagan (See \cite{glen1} and \cite{glen2}). Next, using this result 
together with Cauchon's description for the set of 
$\hc$-invariant prime ideals of $\mnc$ (See \cite{c2}, Th\'eor\`eme 3.2.1), we have constructed 
an algorithm which provides an explicit generating set of quantum minors for each $\hc$-invariant prime ideal in $\mnc$ 
(See \cite{lau2} or \cite{lau3}).

On the other hand, Goodearl and Lenagan have shown (in the general case where $q\in \mathbb{K}^*$ is not a root of unity) 
that, in order to obtain descriptions of all the $\hc$-invariant prime ideals of $R$, we just need to 
determine the $\hc$-invariant prime ideals of certain "localized step-triangular factors" of $R$, namely the algebras
$$R_{\mathbf{r}}^+:= \frac{R}{\langle \yia \mbox{ $\mid$ } \alpha > t \mbox{ or } i < r_{\alpha} \rangle}
\left[ \overline{Y}_{r_1,1}^{-1},\dots,\overline{Y}_{r_t,t}^{-1} \right]$$
and 
$$R_{\mathbf{c}}^-:= \frac{R}{\langle \yia \mbox{ $\mid$ } i > t \mbox{ or } \alpha < c_i \rangle}
\left[ \overline{Y}_{1,c_1}^{-1},\dots,\overline{Y}_{t,c_t}^{-1} \right],$$
where $t \in \gc 0,n \dc$ and where $\mathbf{r}=(r_1,\dots,r_t)$ and $\mathbf{c}=(c_1,\dots,c_t)$ are strictly 
increasing sequences of integers in the range $1,\dots,n$ (See \cite{glen1}, Theorem 3.5). 
Using this result, Goodearl and Lenagan have 
computed the $\hc$-invariant prime ideals of $O_{q}\left( \mathcal{M}_2(\mathbb{K}) \right)$ 
(See \cite{glen1}) and $O_{q}\left( \mathcal{M}_3(\mathbb{K}) \right)$ (See \cite{glen2}).

The aims of this paper are to provide a description for the set $\hc\mbox{-}Spec(R_{\mathbf{r}}^+) $ 
of $\hc$-invariant prime ideals of $R_{\mathbf{r}}^+$ and to count the rank $t$ $\hc$-invariant prime 
ideals of $R$ ($t \in \gc 0,n \dc$), that is 
those $\hc$-invariant prime ideals of $R$ which contain all $(t+1) \times (t+1)$ quantum minors but not all 
$t \times t$ quantum minors. In \cite{glen1}, the authors have shown that the number of rank $t$ $\hc$-invariant prime 
ideals of $R$ is a perfect square. More precisely, they have established (See \cite{glen1}, 3.6) that, for any 
$t \in \gc 0,n \dc$:
\begin{eqnarray}
\label{eqintro}
\mid \hc\mbox{-}Spec^{[t]} (R) \mid \ = 
\left( \sum_{\stackrel{\mathbf{r}=(r_1,\dots,r_t)}{1\leq r_1 < \dots < r_t \leq n }} \mid \hc \mbox{-} Spec(R_{\mathbf{r}}^+) \mid \right)^2
\end{eqnarray}
where $\hc\mbox{-}Spec^{[t]} (R) $ denotes the set of rank $t$ $\hc$-invariant prime ideals of $R$ 
and where $\hc\mbox{-}Spec(R_{\mathbf{r}}^+) $ denotes the set of $\hc$-invariant prime ideals of $R_{\mathbf{r}}^+$. 
The above relation (\ref{eqintro}) opens a potential route to count the rank $t$ $\hc$-invariant prime ideals of $R$: 
if we can compute the number of $\hc$-invariant prime ideals of $R_{\mathbf{r}}^+$, 
then we will be able to count the rank $t$ $\hc$-invariant prime ideals of $R$.

So, to compute the number of rank $t$ $\hc$-invariant prime ideals of $R$, the first step is to study the 
$\hc$-invariant prime ideals of $R_{\mathbf{r}}^+$. Since this algebra is induced from $R$ by factor 
and localization, we first construct (See Section 2), by using the deleting derivations theory (See \cite{c2}),
 $\hc$-invariant prime ideals of $R$ that provide, after factor and localization, 
$2^{r_2-r_1}\dots t^{r_t-r_{t-1}} (t+1)^{n-r_t}$ $\hc$-invariant prime ideals of $R_{\mathbf{r}}^+$ (See Section 
\ref{sectionR+r}). Next, by using (\ref{eqintro}), we are able to show that the number of rank $t$ $\hc$-invariant prime ideals of $R$ 
is greater than or equal to $(t!)^2 S(n+1,t+1)^2$, where $S(n+1,t+1)$ denotes the Stirling number 
of second kind associated to $n+1$ and $t+1$ (See Proposition \ref{propmino}). Finally, after observing that the 
number of $\hc$-invariant prime ideals of $R$ is equal to the poly-Bernoulli number $B_n^{(-n)}$ (See Proposition 
\ref{propnbreHprem}), we use 
a closed formula for the poly-Bernoulli number $B_n^{(-n)}$ (See \cite{kaneko2}, Theorem 2) in order to prove our main result: 
the number of rank $t$ $\hc$-invariant prime ideals of $R$ is actually equal to $(t!)^2 S(n+1,t+1)^2$. 
This result was conjectured by Goodearl, Lenagan and McCammond. As a corollary, we obtain a description 
for the set of $\hc$-invariant prime ideals of $R_{\mathbf{r}}^+$ (See Section \ref{lastsection}). 
\\$ $ 

\section{$\hc$-invariant prime ideals in $\mnk$.}
\label{section1}
$ $

Throughout this paper, we use the following conventions:
\\$ $
\\$\bullet$ If $I$ is a finite set, $|I|$ denotes its cardinality.
\\$\bullet$ $\mathbb{K}$ denotes a (commutative) field and we set $\mathbb{K}^*:=\mathbb{K}\setminus \{0\}$.
\\$\bullet$ \textbf{$q\in \mathbb{K}^*$ is  not a root of unity}.
\\$\bullet$ $n$ denotes a positive integer with $n \geq 2$.
\\$\bullet$ $R=\mnk$ denotes the quantization of the ring of regular functions on 
$n \times n$ matrices with entries in $\mathbb{K}$; it is the $\mathbb{K}$-algebra generated by the 
$n \times n$ indeterminates $Y_{\ia}$, $1 \leq i,\alpha \leq n$, subject to the following 
relations:\\
If $\left( \begin{array}{cc} x & y \\ z & t \end{array} \right)$ is any $2 \times 2$ sub-matrix of 
$\mathcal{Y}:=\left( Y_{\ia} \right) _{(\ia) \in \gc 1,n \dc^2}$, then
 \begin{enumerate}
 \item $ yx=q^{-1} xy, \quad zx=q^{-1} xz, \quad zy=yz, \quad ty=q^{-1}
 yt, \quad tz=q^{-1} zt$.
 \item $tx=xt-(q-q^{-1})yz$.
\end{enumerate}
 These relations agree with the relations used in \cite{c2}, \cite{glen1}, \cite{glen2}, \cite{lau} and \cite{lau2}, 
but they 
differ from those of \cite{pw} and \cite{c3} by an interchange of $q$ and $q^{-1}$. It is well known that $R$ can be 
presented as an iterated Ore extension over $\mathbb{K}$, with the generators $\yia$ adjoined in lexicographic order. 
Thus the ring $R$ is a Noetherian domain. \textbf{We denote by $F$ its skew-field of fractions}. Moreover, since 
$q$ is not a root of unity, it follows from \cite[Theorem 3.2]{gl4} that all prime ideals 
of $R$ are completely prime.
\\$\bullet$ It is well known that the group $\hc:=\left( \mathbb{C}^* \right)^{2n}$ acts 
on $R$ by $\mathbb{K}$-algebra automorphisms via:
$$(a_1,\dots,a_n,b_1,\dots,b_n).Y_{\ia} = a_i b_\alpha Y_{\ia} 
\quad \forall \: (\ia)\in \gc 1,n \dc^2.$$
An \textbf{$\hc $-eigenvector} $x$ of $R$ is a nonzero element $x \in R$ such that $h(x) \in \mathbb{K}^*x$ 
for each $h \in \hc$. An ideal $I$ of $R$ is said to be \textbf{$\hc$-invariant} if $h(I) =I$ for all $h\in \hc$. 
We denote by \underline{$\hc$-$Spec(R)$} the set of $\hc$-invariant prime ideals of $R$. \\$ $

The aim of this paragraph is to construct $\hc$-invariant prime ideals of $R$ that, after factor and localization, will 
provide $\hc$-invariant prime ideals of $R_{\mathbf{r}}^+$ (See the introduction for the definition 
of this algebra). In order to do this, we use the description 
of the set $\hc$-$Spec(R)$ that Cauchon has obtained by applying the theory of deleting derivations (See \cite{c2}).
\\$ $

\subsection{Standard deleting derivations algorithm and description of $\hc$-$Spec(R)$.}
\label{subsection1}
$ $

In this section, we provide the background definitions and notations for the 
standard deleting derivations algorithm (See \cite{c2,lau,lau2}) and we recall the description of 
the set $\hc$-$Spec(R)$ that Cauchon has obtained by using this algorithm (See \cite{c2}). 
\\$ $
\begin{notas}
$ $
\begin{itemize}
\item We denote by $\leq_s$ the lexicographic ordering on $\mathbb{N}^2$. We often call it 
\textbf{the standard ordering on $\mathbb{N}^2$}. Recall that $(\ia) \leq_s (j,\beta) \Longleftrightarrow [(i < j) \mbox{ or } (i=j \mbox{ and } \alpha \leq \beta )]$.
\item We set $E_s=\left(\gc 1,n \dc^2 \cup \{(n,n+1)\} \right) \setminus \{(1,1)\}$.
\item Let $(j,\beta) \in E_s$. If $(j,\beta) \neq (n,n+1)$, $(j,\beta)^{+}$ denotes the smallest element 
(relatively to $\leq_s$) of the set $\left\{ (\ia) \in E_s \mbox{ $\mid$ }(j,\beta) <_s (\ia) \right\}$.
\end{itemize}
\end{notas}
$ $

In \cite{c2}, Cauchon has shown that the theory of deleting derivations (See \cite{c1}) can be 
applied to the iterated Ore extension $R=\mathbb{C}[Y_{1,1}]\dots [Y_{n,n};\sigma_{n,n},\delta_{n,n}]$ 
(where the indices are increasing for $\leq_s$). The corresponding deleting derivations algorithm is called 
\textbf{the standard deleting derivations algorithm}. It consists in the construction, for each $r \in E_s$, of the family 
$(Y_{\ia}^{(r)})_{(\ia) \in \gc 1,n \dc ^2}$ of elements of $F=Fract(R)$, defined as follows: 
\\$ $
\begin{enumerate}
\item \underline{If $r=(n,n+1)$}, then 
$Y_{\ia}^{(n,n+1)}=Y_{\ia}$ for all $(\ia) \in \gc 1,n \dc^2$.\\$ $
\item \underline{Assume that $r=(j,\beta) <_s (n,n+1)$}
and that the $Y_{\ia}^{(r^{+})}$ ($(\ia) \in \gc 1,n \dc ^2$) are already constructed. 
Then, it follows from \cite[Th\'eor\`eme 3.2.1]{c1} that $Y_{j,\beta}^{(r^+)} \neq 0$ and, 
for all $(\ia) \in \gc 1,n \dc ^2$, we have:  
$$Y_{\ia}^{(r)}=\left\{ \begin{array}{ll}
Y_{\ia}^{(r^{+})}-Y_{i,\beta}^{(r^{+})} \left(Y_{j,\beta}^{(r^{+})}\right)^{-1}
 Y_{j,\alpha}^{(r^{+})} 
& \mbox{ if } i<j \mbox{ and } \alpha < \beta \\ 
Y_{\ia}^{(r^{+})} & \mbox{ otherwise.}
\end{array} \right.$$
\end{enumerate}

\begin{nota}
$ $
\\Let $r \in E_s$. We denote by $R^{(r)}$ the subalgebra of $F=Fract(R)$ generated by the $Y_{\ia}^{(r)}$ 
($(\ia) \in \gc 1,n \dc^2$), that is, $R^{(r)}:=\mathbb{C} \langle Y_{\ia}^{(r)} \mbox{ $\mid$ }
(\ia) \in \gc 1,n \dc^2  \rangle$. 
\end{nota}
$ $

\begin{notas}
$ $
\\We set $\ov{R}:=R^{(1,2)}$ and 
$T_{\ia}:=Y_{\ia}^{(1,2)}$ for all $(\ia) \in \gc 1,n \dc^2$.
\end{notas}
$ $

Let $(j,\beta) \in E_s$ with $(j,\beta) \neq (n,n+1)$. The theory of deleting derivations 
allows us to construct embeddings $\varphi_{(j,\beta)}:Spec(R^{(j,\beta)^{+}}) \longrightarrow Spec(R^{(j,\beta)})$
 (See \cite{c1}, 4.3). By composition, we obtain an embedding 
$\varphi:Spec(R) \longrightarrow Spec(\ov{R})$ which is called \textbf{the canonical embedding}. In \cite{c2}, 
Cauchon has described the set $\hc$-$Spec(R)$ by determining its "canonical image" $\varphi(\hc$-$Spec(R))$. To do this, 
he has introduced the following conventions and notations.
\\$ $
\begin{convs}
\label{convstronq}
$ $
\begin{itemize}
\item Let $v=(l,\gamma) \in \gc 1,n \dc^2$.
\begin{enumerate}
\item The set $C_v:=\{(i,\gamma) \mbox{ $\mid$ } 1 \leq i \leq l \} \subset \gc 1,n \dc^2$ is called 
the \textbf{truncated column with extremity $v$}.
\item The set $L_v:=\{(l,\alpha) \mbox{ $\mid$ } 1 \leq \alpha \leq \gamma \} \subset \gc 1,n \dc^2$ is called 
the \textbf{truncated row with extremity $v$}.
\end{enumerate}
\item $W$ denotes the set of all the subsets in $\gc 1,n \dc^2 $ which are a union of truncated rows and columns. 
\end{itemize}
\end{convs}
$ $

\begin{nota}
$ $
\\Given $w \in W$, $K_w$ denotes the ideal in $\ov{R}$ generated by the $T_{\ia}$ such that $(\ia) \in w$.
\\(Recall that $K_w$ is a completely prime ideal in the quantum affine space $\ov{R}$ (See \cite{gl2}, 2.1).) 
\end{nota}
 $ $

The following description of the set $\hc$-$Spec(R)$ was obtained by Cauchon (See \cite{c2}, 
Corollaire 3.2.1). 
\\$ $
\begin{prop}
\label{descriptionHprem}
$ $
\begin{enumerate}
\item Given $w \in W$, there exists a (unique) $\hc$-invariant (completely) prime ideal $J_w$ in $R$ such 
that $\varphi(J_w)=K_w$.
\item $ \hc \mbox{-}Spec(R) = \{J_w \mbox{ $\mid$ }w\in W\}$. 
\end{enumerate}
\end{prop}
$ $

\subsection{Number of $\hc$-invariant prime ideals in $R$.}
$ $

In \cite{c2}, Cauchon has used his description of the set $\hc$-$Spec(R)$ in order to 
give a formula for the total number $S(n)$ of $\hc$-invariant prime ideals of $R$. More precisely, 
he has established (See \cite{c2}, Proposition 3.3.2) that:
$$ S(n)=(-1)^{n-1} \sum_{k=1}^n (k+1)^n \sum_{j=1}^k  (-1)^{j-1} 
\left( \begin{array}{l} k \\ j \end{array} \right) j^n,
$$
that is
$$S(n)=(-1)^{n} \sum_{k=1}^n (-1)^k k! (k+1)^n \left( \frac{(-1)^k}{k!}\sum_{j=1}^k  (-1)^{j} 
\left( \begin{array}{l} k \\ j \end{array} \right) j^n \right).
$$ 
Recall (See \cite{stanleybook}, p. 34) that $\displaystyle{\frac{(-1)^k}{k!}\sum_{j=1}^k  (-1)^{j} 
\left( \begin{array}{l} k \\ j \end{array} \right) j^n=\frac{(-1)^k}{k!}\sum_{j=0}^k  (-1)^{j} 
\left( \begin{array}{l} k \\ j \end{array} \right) j^n}$ is equal to the Stirling number of second kind
 $S(n,k)$ (See, for example, \cite{stanleybook} for more details on the Stirling numbers of second kind). 
Hence, we have: 
$$S(n)=(-1)^{n} \sum_{k=1}^n (-1)^k k! (k+1)^n S(n,k),$$
that is
\begin{eqnarray}
\label{eqtotal}
S(n)=(-1)^{n} \sum_{k=1}^n \frac{(-1)^k k!}{(k+1)^{-n} }S(n,k).
\end{eqnarray}
On the other hand, it follows from \cite[Theorem 1]{kaneko1} that: 
$$(-1)^{n} \sum_{k=0}^n \frac{(-1)^k k!}{(k+1)^{-n} }S(n,k)=B_n^{(-n)},$$
where $B_n^{(-n)}$ denotes the poly-Bernoulli number associated to $n$ and $-n$ (See \cite{kaneko1} for the definition 
of the poly-Bernoulli numbers). Observing that $S(n,0)=0$ (See \cite{stanleybook}), we get:
$$(-1)^{n} \sum_{k=1}^n \frac{(-1)^k k!}{(k+1)^{-n} }S(n,k)=B_n^{(-n)},$$
and thus, we deduce from (\ref{eqtotal}) that: 
\\$ $

\begin{prop}
\label{propnbreHprem}
$$\mid \hc \mbox{-}Spec(R) \mid \ =B_n^{(-n)}.$$
\end{prop}
$ $

This rewriting of Cauchon's formula was first obtained by Goodearl and McCammond.
\\$ $

\subsection{Vanishing and non-vanishing criteria for the entries of $q$-quantum matrices.} 
\label{section2}
$ $

Let $J_w$ ($w \in W$) be an $\hc$-invariant prime ideal of $R$ (See Proposition \ref{descriptionHprem}). 
In the next section, we will need to know which indeterminates $\yia$ belong to $J_w$, 
that is which $y_{\ia}:=\yia+J_w$ are zero. This problem is dealt with in 
Proposition \ref{propyiapasdansJw} and Proposition \ref{propannulationtiawr} where we respectively 
obtain a non-vanishing criterion and a vanishing criterion for the entries of $q$-quantum 
matrices.
\\$ $

For the remainder of this section, $K$ denotes a $\mathbb{K}$-algebra which is also a skew-field. 
Except otherwise stated, all the considered matrices have their entries in $K$.
\\$ $

\begin{defis}
$ $
\\Let $M=(x_{i,\alpha})_{(\ia) \in \gc 1,n \dc^2}$
 be a $n \times n$ matrix and let $(j,\beta) \in E_s$.
\begin{itemize}
\item We say that $M$ is a \textbf{$q$-quantum matrix} if the following relations hold between 
the entries of $M$:
 \\If $\left( \begin{array}{cc} x & y \\ z & t \end{array} \right)$ is any $2 \times 2$
 sub-matrix of $M$, then 
 \begin{enumerate}
 \item $ yx=q^{-1} xy, \quad zx=q^{-1} xz, \quad zy=yz, \quad ty=q^{-1}
 yt, \quad tz=q^{-1} zt.$
 \item $tx=xt-(q-q^{-1})yz$.
\end{enumerate}

\item We say that $M$ is a \textbf{$(j,\beta)$-$q$-quantum matrix} 
if the following relations hold between the entries of $M$:
 \\If $\left( \begin{array}{cc} x & y \\ z & t \end{array} \right)$ is any $2 \times 2$
 sub-matrix of $M$, then 
 \begin{enumerate}
 \item $ yx=q^{-1} xy, \quad zx=q^{-1} xz, \quad zy=yz, \quad ty=q^{-1}
 yt, \quad tz=q^{-1} zt.$
 \item If $t=x_v$, then $\left\{ \begin{array}{lll}
v \geq_s (j,\beta) & \Longrightarrow & tx=xt \\
v <_s (j,\beta) & \Longrightarrow & tx=xt-(q-q^{-1})yz. \\
\end{array} \right.$
\end{enumerate}
 \end{itemize}
\end{defis}
$ $

\begin{convs}
\label{conv1}
$ $
\\Let $M=(x_{\ia})_{(\ia) \in \gc 1,n \dc^2}$ be a $q$-quantum matrix. 
\\As $r$ runs over the set $E_s$, we define matrices $M^{(r)} =(x_{\ia}^{(r)})_{(\ia) \in \gc 1,n \dc^2}$ 
as follows:
\begin{enumerate}
\item \underline{If $r=(n,n+1)$}, then the entries of the matrix $M^{(n,n+1)}$ are defined by 
$x_{\ia}^{(n,n+1)}:=x_{\ia}$ for all $(\ia) \in \gc 1,n \dc^2$. 
\item \underline{Assume that $r=(j,\beta) \in E_s \setminus \{(n,n+1)\}$} and that the matrix $M^{(r^{+})}$ is 
already known. The entries $x_{\ia}^{(r)}$ of the matrix $M^{(r)}$ are defined as follows:
\begin{enumerate}
\item If $x_{j,\beta}^{(r^+)}=0$, then $x_{\ia}^{(r)}=x_{\ia}^{(r^+)}$ 
for all $(\ia) \in \gc 1,n \dc^2$.
\item If $x_{j,\beta}^{(r^+)}\neq 0$ and $(\ia) \in \gc 1,n \dc^2$, then 
\\$x_{\ia}^{(r)}= \left\{ \begin{array}{ll}
x_{\ia}^{(r^+)}-x_{i,\beta}^{(r^+)} \left( x_{j,\beta}^{(r^+)}\right)^{-1} 
x_{j,\alpha}^{(r^+)}
& \qquad \mbox{if }i <j \mbox{ and } \alpha < \beta \\
x_{\ia}^{(r^+)} & \qquad \mbox{otherwise.} \end{array} \right.$ 
\end{enumerate}
We say that \textbf{$M^{(r)}$ is the matrix obtained from $M$ by applying the standard deleting derivations 
algorithm at step $r$}.
\item \underline{If $r=(1,2)$}, we set $t_{\ia}:=x_{\ia}^{(1,2)}$ for all $(\ia) \in \gc 1,n \dc^2$.
\end{enumerate}
\end{convs}
$ $

Observe that the formulas of Conventions \ref{conv1} allow us to 
express the entries of $M^{(r^+)}$ in terms of those of $M^{(r)}$.
\\$ $

\begin{prop}[Restoration algorithm]
\label{algorestitution}
$ $
\\Let $M=(x_{\ia})_{(\ia) \in \gc 1,n \dc^2 }$ be a $q$-quantum matrix and let 
$r=(j,\beta) \in E_s$ with $r \neq (n,n+1)$. 
\begin{enumerate}
\item If $x_{j,\beta}^{(r)}=0$, then $x_{\ia}^{(r^+)}=x_{\ia}^{(r)}$ 
for all $(\ia) \in \gc 1,n \dc^2$.
\item If $x_{j,\beta}^{(r)}\neq 0$ and $(\ia) \in \gc 1,n \dc^2$, then 
\\$x_{\ia}^{(r^+)}= \left\{ \begin{array}{ll}
x_{\ia}^{(r)}+x_{i,\beta}^{(r)} \left( x_{j,\beta}^{(r)}\right)^{-1} 
x_{j,\alpha}^{(r)}
& \qquad \mbox{if }i <j \mbox{ and } \alpha < \beta \\
x_{\ia}^{(r)} & \qquad \mbox{otherwise.} \end{array} \right.$ 
\end{enumerate}
\end{prop}
$ $

Note that our definitions of $q$-quantum matrix and $(j,\beta)$-$q$-quantum matrix slightly 
differ from those of \cite{c3} (See \cite{c3}, D\'efinitions III.1.1 and III.1.3). Because of this, 
we must interchange $q$ and $q^{-1}$ whenever carrying over result of \cite{c3}. 
\\$ $

\begin{lem}
\label{lemjbetasqquantique}
$ $
\\Let $(j,\beta) \in E_s$.
\\If $M=(x_{i,\alpha})_{(\ia) \in \gc 1,n \dc^2}$ is a $q$-quantum matrix, then the matrix $M^{(j,\beta)}$ 
is $(j,\beta)$-$q$-quantum.
\end{lem}
\preuve This lemma is proved in the same manner as \cite[Proposition III.2.3.1]{c3}. \fin
\\$ $

We deduce from the above Lemma \ref{lemjbetasqquantique} the following non-vanishing criterion for the entries of 
a $q$-quantum matrix.
\\$ $

\begin{prop}
\label{propyiapasdansJw}
$ $
\\Let $M=(x_{i,\alpha})_{(\ia) \in \gc 1,n \dc^2}$ be a $q$-quantum matrix and let $(\ia) \in \gc 1,n \dc^2$.
\\If $t_{\ia} \neq 0$, then $x_{\ia} \neq 0$. In other words, if $x_{\ia} =0$, then $t_{\ia}=0$.
\end{prop}
\preuve Assume that $x_{\ia} =0$. We first prove that $x_{\ia}^{(j,\beta)}=0$ for all $(j,\beta) \in E_s$. 
To achieve this aim, we proceed by decreasing induction (for $\leq_s$) on $(j,\beta)$.

Since $x_{\ia}^{(n,n+1)}=x_{\ia}$, the case $(j,\beta)=(n,n+1)$ is done. Assume now that 
$(j,\beta) <_s (n,n+1)$ and $x_{\ia}^{(j,\beta)^+}=0$. If $x_{\ia}^{(j,\beta)}=x_{\ia}^{(j,\beta)^+}$, 
we obviously have $x_{\ia}^{(j,\beta)}=0$. Next, if $x_{\ia}^{(j,\beta)} \neq x_{\ia}^{(j,\beta)^+}$, 
then $i < j$ and $\alpha < \beta$. Hence, it follows from Lemma \ref{lemjbetasqquantique} that 
the matrix $\displaystyle{\left( \begin{array}{cc} 
x_{\ia}^{(j,\beta)^+} & x_{i,\beta}^{(j,\beta)^+} \\
x_{j,\alpha}^{(j,\beta)^+} & x_{j,\beta}^{(j,\beta)^+}
\end{array} \right)}$ is $q$-quantum, so that 
$$x_{j,\beta}^{(j,\beta)^+}x_{\ia}^{(j,\beta)^+} - 
x_{\ia}^{(j,\beta)^+} x_{j,\beta}^{(j,\beta)^+}=-(q-q^{-1}) 
x_{i,\beta}^{(j,\beta)^+}x_{j,\alpha}^{(j,\beta)^+}.$$
Since $x_{\ia}^{(j,\beta)^+}=0$, we deduce from this equality that, in $K$, 
$x_{i,\beta}^{(j,\beta)^+}x_{j,\alpha}^{(j,\beta)^+}=0$. Thus, 
$x_{i,\beta}^{(j,\beta)^+} = 0$ or $x_{j,\alpha}^{(j,\beta)^+}=0$. On the other hand, since $i<j$ and $\alpha < \beta$, we have 
$x_{\ia}^{(j,\beta)}=x_{\ia}^{(j,\beta)^+}-x_{i,\beta}^{(j,\beta)^+} 
\left( x_{j,\beta}^{(j,\beta)^+}\right)^{-1} x_{j,\alpha}^{(j,\beta)^+}$. 
Now it follows from the induction hypothesis that $x_{\ia}^{(j,\beta)^+}=0$. Hence, we have 
\\$x_{\ia}^{(j,\beta)}= -x_{i,\beta}^{(j,\beta)^+} \left( x_{j,\beta}^{(j,\beta)^+}\right)^{-1} x_{j,\alpha}^{(j,\beta)^+}$. Finally, since 
$x_{i,\beta}^{(j,\beta)^+} = 0$ or $x_{j,\alpha}^{(j,\beta)^+}=0$, we get $x_{\ia}^{(j,\beta)}=0$, 
as desired. This achieves the induction.

In particular, we have shown that $x_{\ia}^{(1,2)}=0$, that is $t_{\ia}=0$. \fin
\\$ $

Proposition \ref{propyiapasdansJw} furnishes a non-vanishing criterion for the entries of a 
$q$-quantum matrix. In order to construct, in the next section, $\hc$-invariant prime ideals of $R$ that 
will provide, after factor and localization, $\hc$-invariant prime ideals of 
$\displaystyle{R_{\mathbf{r}}^+:= \frac{R}{\langle \yia \mbox{ $\mid$ } \alpha > t \mbox{ or } i < r_{\alpha} \rangle}
\left[ \overline{Y}_{r_1,1}^{-1},\dots,\overline{Y}_{r_t,t}^{-1} \right]}$ ($\mathbf{r}=(r_1,\dots,r_t)$ 
with $1 \leq r_1 < \dots < r_t \leq n $), we also need to get a vanishing criterion 
for the entries $x_{\ia}$, $\alpha > t \mbox{ or } i < r_{\alpha}$, of a $q$-quantum matrix. This is what we do now.
\\$ $

\begin{nota}
$ $
\\If $t$ denotes an element of $\gc 0,n \dc$, we set:
$$\mathbf{R}_t:=\{(r_1,\dots,r_t) \in \mathbb{N} \mbox{ $\mid$ } 1 \leq r_1 < \dots < r_t \leq n \}.$$
(If $t=0$, then $\mathbf{R}_0= \emptyset$.) 
\end{nota}
$ $

For the remainder of this section, we fix $t \in \gc 0,n \dc$ and $\mathbf{r}=(r_1,\dots,r_t) \in \mathbf{R}_t$, 
and we denote by $w_{\mathbf{r}}$ the subset of $\gc 1,n \dc ^2$ corresponding to indeterminates $\yia$ that have been 
set equal to zero in $R_{\mathbf{r}}^+$, that is, we set: 
$$w_{\mathbf{r}}:=\left[ \bigcup_{\alpha \in \gc 1,t \dc} \gc 1,r_{\alpha} -1 \dc \times \{ \alpha \} 
\right] \bigcup \gc 1,n \dc \times \gc t+1,n \dc .$$
$ $

For instance, if $n=3$, $t=2$ and $\mathbf{r}=(1,3)$, we have: 
\\$ $
\begin{center}
$w_{(1,3)}=\begin{tabular}{| b{1mm} | b{1mm} | b{1mm} |}
    \hline
     \LCC &  \black &  \black \\ 
           & &  \ECC \\    
    \hline
     \LCC  & \black & \black  \\
      & &  \ECC \\
    \hline
      \LCC & & \black \\
      & &  \ECC \\
    \hline
      \end{tabular} $
, where the black boxes symbolize the elements of $w_{(1,3)}$.
\end{center}
$ $

Note that $w_{\mathbf{r}}$ is a union of truncated columns, so that:
\\$ $
\begin{rem}
$ $
\\$w_{\mathbf{r}}$ belongs to $W$.  
\end{rem}
$ $


\begin{obs}
\label{obswrconvexe}
$ $
\\Let $(\ia) \in w_{\mathbf{r}}$. If $\beta \in \gc \alpha ,n \dc$, then $(i,\beta) \in w_{\mathbf{r}}$.
\end{obs}
\preuve We distinguish two cases.
\\$\bullet$ If $(\ia) \in \gc 1,n \dc \times \gc t+1,n \dc$, then $\alpha \geq t+1$. Hence $\beta \geq \alpha \geq t+1$ 
and thus, we have $(i,\beta) \in \gc 1,n \dc \times \gc t+1,n \dc \subseteq w_{\mathbf{r}}$, as required.
\\$\bullet$ Assume now that 
$\displaystyle{(\ia) \in \bigcup_{\gamma \in \gc 1,t \dc} \gc 1,r_{\gamma} -1 \dc \times \{ \gamma \}}$, so that 
we have $\alpha \leq t$ and $i \leq r_{\alpha} -1$. If $\beta > t$, we conclude
 as in the previous case that $(i,\beta) \in w_{\mathbf{r}}$. So we assume that $\beta \leq t$. Since 
$i \leq r_{\alpha} -1$ and since $\alpha \leq \beta \leq t$, we have $i \leq r_{\alpha } -1 \leq r_{\beta} -1$. 
Hence, $(i,\beta) \in \gc 1, r_{\beta}-1 \dc \times \{\beta\} \subseteq w_{\mathbf{r}}$, as desired. \fin
\\$ $

This observation allows us to prove the following vanishing criterion:
\\$ $

\begin{prop}
\label{propannulationtiawr}
$ $
\\Let $M=(x_{\ia})_{(\ia) \in \gc 1,n \dc^2}$ be a $q$-quantum matrix. 
\\If $t_{\ia} =0$ for all $(\ia) \in w_{\mathbf{r}}$, then  
$x_{\ia} =0$ for all $(\ia) \in w_{\mathbf{r}}$.
\end{prop}  
\preuve Assume that $t_{\ia} =0$ for all $(\ia) \in w_{\mathbf{r}}$.
 We first prove by induction on $(j,\beta)$ (with respect of $\leq_s$) that $x_{\ia}^{(j,\beta)} =0$
 for all $(\ia) \in w_{\mathbf{r}}$ and $(j,\beta) \in E_s$.

If $(j,\beta)=(1,2)$, then $x_{\ia}^{(1,2)}=t_{\ia}=0$ for all $(\ia) \in w_{\mathbf{r}}$, as required. 
Assume now that $(j,\beta) <_s (n,n+1)$ and that $x_{\ia}^{(j,\beta)} =0$ for all $(\ia) \in w_{\mathbf{r}}$. 
Let $(\ia) \in w_{\mathbf{r}}$. If $x_{\ia}^{(j,\beta)^+}=x_{\ia}^{(j,\beta) }$, the desired result follows from 
the induction hypothesis. Next, if $x_{\ia}^{(j,\beta)^+} \neq x_{\ia}^{(j,\beta) }$, 
it follows from Proposition \ref{algorestitution} that 
$x_{j,\beta}^{(j,\beta)} \neq 0$, $i <j$, $\alpha < \beta$ and 
$x_{\ia}^{(j,\beta)^+}= x_{\ia}^{(j,\beta)}+x_{i,\beta}^{(j,\beta)} \left( x_{j,\beta}^{(j,\beta)}\right)^{-1} 
x_{j,\alpha}^{(j,\beta)}$. Since $(\ia) \in w_{\mathbf{r}}$, we deduce from the induction hypothesis 
that $x_{\ia}^{(j,\beta)} =0$, so that $x_{\ia}^{(j,\beta)^+}= 
x_{i,\beta}^{(j,\beta)} \left( x_{j,\beta}^{(j,\beta)}\right)^{-1} x_{j,\alpha}^{(j,\beta)}$. 
Moreover, since $(\ia) \in w_{\mathbf{r}}$ and $\alpha < \beta$, it follows from Observation \ref{obswrconvexe} 
that $(i,\beta) \in w_{\mathbf{r}}$. Then, we deduce from the induction hypothesis that 
$x_{i,\beta}^{(j,\beta)} =0$, so that $x_{i,\alpha}^{(j,\beta)^+}=x_{i,\beta}^{(j,\beta)} 
\left( x_{j,\beta}^{(j,\beta)}\right)^{-1} x_{j,\alpha}^{(j,\beta)}=0$. This achieves the induction.

In particular, we have proved that $x_{\ia}=x_{\ia}^{(n,n+1)}=0$ for all $(\ia) \in  w_{\mathbf{r}}$. \fin
\\$ $ 

\subsection{$\hc$-invariant prime ideals $J_w$ with $w_{\mathbf{r}} \subseteq w$.}
\label{subsection24}
$ $

As in the previous section, we fix $t \in \gc 0,n \dc$ and $\mathbf{r}=(r_1,\dots,r_t) \in \mathbf{R}_t$, and 
we set: 
$$w_{\mathbf{r}}:= 
\left[ \bigcup_{\alpha \in \gc 1,t \dc} \gc 1,r_{\alpha} -1 \dc \times \{ \alpha \} \right] 
\bigcup \gc 1,n \dc \times \gc t+1,n \dc  .$$
Recall (See Proposition \ref{descriptionHprem}) that, if $w \in W$, 
there exists a (unique) $\hc$-invariant prime ideal of $R$ associated to $w$ 
(See Proposition \ref{descriptionHprem}) and that the $J_w$ 
($w \in W$) are exactly the $\hc$-invariant prime ideals in $R$. 
This section is devoted to the $\hc$-invariant prime ideals $J_w$ ($w\in W$) of $R$ with $w_{\mathbf{r}} \subseteq w$. 
More precisely, we want to know which indeterminates $\yia$ belong to these ideals.

\begin{notas}
\label{notasmw}
$ $
\\Let $w\in W$.
\begin{enumerate}
\item Set $\displaystyle{R_w:=\frac{R}{J_w}}$. It follows from \cite[Lemme 5.3.3]{c1} that, using
 the notations of Section \ref{subsection1}, $R_w$ 
and $\displaystyle{\frac{\ov{R}}{K_w}}$ are two Noetherian algebras with no zero-divisors, 
which have the same skew-field of fractions. We set $F_w:=Fract\left(R_w \right)=Fract\left( \displaystyle{\frac{\ov{R}}{K_w}}  \right)$.
\item If $(\ia) \in \gc 1,n \dc^2$, $y_{\ia}$ denotes the element of $R_w$ defined by 
$y_{\ia}:=\yia+J_w$.
\item We denote by $M_w$ the matrix, with entries in the $\mathbb{K}$-algebra $F_w$, defined by:
$$M_w:=\left( y_{\ia} \right)_{(\ia) \in \gc 1,n \dc^2}.$$
\end{enumerate}
\end{notas}
$ $

Let $w \in W$. Since $\mathcal{Y}=\left( Y_{\ia} \right)_{(\ia) \in \gc 1,n \dc^2}$ 
is a $q$-quantum matrix, $M_w$ is also a $q$-quantum matrix. 
Thus, we can apply the standard deleting derivations algorithm to $M_w$ (See Conventions \ref{conv1} 
with $K=F_w$) and if we still denote $t_{\ia}:=y_{\ia}^{(1,2)}$ for $(\ia) \in \gc 1,n \dc^2$, we get: 
\\$ $

\begin{prop}
\label{theo1}
$ $
\\$t_{\ia}=0$ if and only if $(\ia) \in w$. 
\end{prop}
\preuve By \cite[Propositions 5.4.1 and 5.4.2]{c1}, there exists a $\mathbb{K}$-algebra homomorphism
 $f_{(1,2)}:\ov{R} \rightarrow F_w$ 
such that $f_{(1,2)} \left( T_{\ia} \right) =t_{\ia}$ for $(\ia) \in \gc 1,n \dc^2$. Its kernel is $K_w$ 
and its image is the subalgebra of $F_w$ generated by the $t_{\ia}$ with $(\ia) \in \gc 1,n \dc^2$. 
Hence, $t_{\ia}=0$ if and only if $T_{\ia}\in K_w$, that is, if and only if $(\ia) \in w$. \fin
\\$ $

Consider now an element $w$ in $W$ with $w_{\mathbf{r}} \subseteq w$ and denote by $J_w$ 
the (unique) $\hc$-invariant prime ideal of $R$ associated to $w$ (See Proposition \ref{descriptionHprem}).
 Since $w_{\mathbf{r}} \subseteq w$, we deduce from Proposition \ref{theo1} that $t_{\ia} =0$ 
for all $(\ia) \in w_{\mathbf{r}}$. Hence, we can apply Proposition \ref{propannulationtiawr} to 
the $q$-quantum matrix $M_w$ and we obtain that $y_{\ia}=0$ for all $(\ia) \in w_{\mathbf{r}}$, that is, 
$Y_{\ia}\in J_w$ for all $(\ia) \in w_{\mathbf{r}}$. So we have just established: 
\\$ $

\begin{prop}
\label{propyiainJw}
$ $
\\Let $w \in W$ with $w_{\mathbf{r}} \subseteq w$. If $(\ia) \in w_{\mathbf{r}}$, then $Y_{\ia} $ belongs to $J_w$.
\end{prop}
$ $

We will now add truncated rows to the "$w_{\mathbf{r}}$ diagram" in order to obtain $\hc$-invariant prime 
ideals of $R$ that will provide, after factor and localisation, $\hc$-invariant prime ideals 
of $R_{\mathbf{r}}^+$. We will see later (See Section \ref{lastsection}) that the 
$\hc$-invariant prime ideals of $R$ obtained by adding truncated rows to the "$w_{\mathbf{r}}$ diagram" 
are the only $\hc$-invariant prime ideals of $R$ that will provide, 
after factor and localisation, $\hc$-invariant prime ideals 
of $R_{\mathbf{r}}^+$.
\\$ $

\begin{nota}
\label{notagammar}
$ $
\\We set $\Gamma_{\mathbf{r}}:=\{(\gamma_1,\dots,\gamma_n) \in \mathbb{N}^n \mbox{ $\mid$ } 
\gamma_k \in \gc 0, l \dc \mbox{ if } k \in \gc r_l+1,r_{l+1} \dc\}$.
(Here $r_0=0$ and $r_{t+1}=n$.)
\end{nota}
$ $

For instance, if $n=3$, $t=2$ and $\mathbf{r}=(1,3)$, we have: 
$$\Gamma_{\mathbf{r}}=\{(\gamma_1,\gamma_2,\gamma_3)\in \mathbb{N}^3 \mbox{ $\mid$ } \gamma_1=0, 
\mbox{ } \gamma_2 \leq 1 \mbox{ and } \gamma_3 \leq 1\}.$$
\\$ $

\begin{theo}
\label{theofond}
$ $
\\Let $(\gamma_1,\dots,\gamma_n) \in \Gamma_{\mathbf{r}}$ and 
set $\displaystyle{w_{\mathbf{r},(\gamma_1,\dots,\gamma_n)}:=
w_{\mathbf{r} } \bigcup \left( \bigcup_{k \in \gc 1,n \dc} \{ k \} \times \gc 1,\gamma_k \dc \right) }$.
\\Then $w_{\mathbf{r},(\gamma_1,\dots,\gamma_n)}$ belongs to $W$ and the $\hc$-invariant prime ideal 
$J_{w_{\mathbf{r},(\gamma_1,\dots,\gamma_n)}}$ of $R$ has the following properties:
\begin{enumerate}
\item $\yia \in J_{w_{\mathbf{r},(\gamma_1,\dots,\gamma_n)}}$ for all $(\ia) \in w_{\mathbf{r}}$.
\item $Y_{r_k,k} \notin J_{w_{\mathbf{r},(\gamma_1,\dots,\gamma_n)}}$ for all $k \in \gc 1,t \dc$.
\end{enumerate}
\end{theo} 
\preuve Since $w_{\mathbf{r}}$ is a union of truncated columns and since 
$\displaystyle{\bigcup_{k \in \gc 1,n \dc} \{k\} \times \gc 1,\gamma_k \dc}$ is a union of truncated rows, 
$w_{\mathbf{r},(\gamma_1,\dots,\gamma_n)}$ is a union of truncated rows and columns, so that 
$w_{\mathbf{r},(\gamma_1,\dots,\gamma_n)} \in W$.
\\$ $

Since $w_{\mathbf{r}} \subseteq w_{\mathbf{r},(\gamma_1,\dots,\gamma_n)}$, we deduce from Proposition 
\ref{propyiainJw} that $\yia \in J_{w_{\mathbf{r},(\gamma_1,\dots,\gamma_n)}}$ for all $(\ia) \in w_{\mathbf{r}}$.
\\$ $

Now we want to prove that 
$Y_{r_k,k} \notin J_{w_{\mathbf{r},(\gamma_1,\dots,\gamma_n)}}$ for all $k \in \gc 1,t \dc$. Assume this is not the case, 
that is, assume that there exists $k \in \gc 1,t \dc$ with 
$Y_{r_k,k} \in J_{w_{\mathbf{r},(\gamma_1,\dots,\gamma_n)}}$. Then, $y_{r_k,k}=0$ and it follows from Proposition 
\ref{propyiapasdansJw} that $y_{r_k,k}^{(1,2) }=t_{r_k,k}=0$. Thus, we deduce from 
Proposition \ref{theo1} that $(r_k,k) \in w_{\mathbf{r},(\gamma_1,\dots,\gamma_n)}$.

Observe now that, since $k \leq t$, $(r_k,k) \notin 
\gc 1,n \dc \times \gc t+1,n \dc$. Further, it is obvious that 
$\displaystyle{(r_k,k) \notin \bigcup_{\alpha \in \gc 1,t \dc} \gc 1,r_{\alpha} -1 \dc \times \{ \alpha \}}$. Hence, 
$(r_k,k) \notin w_{\mathbf{r}}$.

All this together shows that $\displaystyle{(r_k,k) \in w_{\mathbf{r},(\gamma_1,\dots,\gamma_n)} \setminus 
w_{\mathbf{r}}=\bigcup_{l \in \gc 1,n \dc} \{ l \} \times \gc 1,\gamma_l \dc }$, so that 
$k \leq \gamma_{r_k}$.

However, since $(\gamma_1,\dots,\gamma_n) \in \Gamma_{\mathbf{r}}$, we have
 $\gamma_{r_k} \leq k-1$. This is a contradiction and thus we have proved that 
$Y_{r_k,k} \notin J_{w_{\mathbf{r},(\gamma_1,\dots,\gamma_n)}}$ for all $k \in \gc 1,t \dc$.  \fin
\\$ $

Let us now give an example for the elements 
$w_{\mathbf{r},(\gamma_1,\dots,\gamma_n)}$ ($(\gamma_1,\gamma_2,\gamma_3) \in \Gamma_{\mathbf{r}}$) of 
Theorem \ref{theofond}.
\\If $n=3$, $t=2$ and $\mathbf{r}=(1,3)$, we have already 
note that $$\Gamma_{\mathbf{r}}=\{(\gamma_1,\gamma_2,\gamma_3)\in \mathbb{N}^3 \mbox{ $\mid$ } \gamma_1=0, 
\mbox{ } \gamma_2 \leq 1 \mbox{ and } \gamma_3 \leq 1\},$$
so that the elements 
$w_{\mathbf{r},(\gamma_1,\dots,\gamma_n)}$ ($(\gamma_1,\gamma_2,\gamma_3) \in \Gamma_{\mathbf{r}}$) of 
Theorem \ref{theofond} are: 
\\$ $
\begin{center}
\begin{tabular}{rr}
$w_{(1,3),(0,0,0)}=w_{(1,3)}=\begin{tabular}{| b{1mm} | b{1mm} | b{1mm} |}
    \hline
     \LCC &  \black &  \black \\ 
           & &  \ECC \\    
    \hline
     \LCC  & \black & \black  \\
      & &  \ECC \\
    \hline
      \LCC & & \black \\
      & &  \ECC \\
    \hline
      \end{tabular} $ &
$w_{(1,3),(0,1,0)}=\begin{tabular}{| b{1mm} | b{1mm} | b{1mm} |}
    \hline
     \LCC &  \black &  \black \\ 
           & &  \ECC \\    
    \hline
     \LCC  \black & \black & \black  \\
      & &  \ECC \\
    \hline
      \LCC & & \black \\
      & &  \ECC \\
    \hline
      \end{tabular} $ \\
 & \\
$w_{(1,3),(0,0,1)}=\begin{tabular}{| b{1mm} | b{1mm} | b{1mm} |}
    \hline
     \LCC &  \black &  \black \\ 
           & &  \ECC \\    
    \hline
     \LCC  & \black & \black  \\
      & &  \ECC \\
    \hline
      \LCC \black & & \black \\
      & &  \ECC \\
    \hline
      \end{tabular} $ &
$w_{(1,3),(0,1,1)}=\begin{tabular}{| b{1mm} | b{1mm} | b{1mm} |}
    \hline
     \LCC &  \black &  \black \\ 
           & &  \ECC \\    
    \hline
     \LCC \black & \black & \black  \\
      & &  \ECC \\
    \hline
      \LCC \black & & \black \\
      & &  \ECC \\
    \hline
      \end{tabular} $
\end{tabular}
\end{center}
$ $
\\(As previously, if $w \in W$, the black boxes symbolize the elements of $w$.)
\\$ $

\section{Number of rank $t$ $\hc$-invariant prime ideals in $\mnk$.}
$ $

In this paragraph, using the previous section, we begin by constructing 
$\hc$-invariant prime ideals of the algebra 
$\displaystyle{R_{\mathbf{r}}^+:= \frac{\mnk}{\langle \yia \mbox{ $\mid$ } \alpha > t \mbox{ or } i < r_{\alpha} \rangle}
\left[ \overline{Y}_{r_1,1}^{-1},\dots,\overline{Y}_{r_t,t}^{-1} \right]}$, 
where $t \in \gc 0,n \dc$ and $\mathbf{r} =(r_1,\dots,r_t)$ is a strictly increasing sequence of integers in the range 
$1,\dots,n$. Next, following the route sketched in the introduction, we establish our main result:
 the number $\mid \hc\mbox{-}Spec^{[t]} (R) \mid $ 
of $\hc$-invariant prime ideals of $R=\mnk$ which contain all $(t+1) \times (t+1)$ quantum minors but not all 
$t \times t$ quantum minors is equal to $(t!)^2 S(n+1,t+1)^2$, where $S(n+1,t+1)$ denotes the Stirling number 
of second kind associated to $n+1$ and $t+1$. From this result, we derive a description 
of the set of $\hc$-invariant prime ideals of $R_{\mathbf{r}}^+$.
\\$ $

\subsection{$\hc$-invariant prime ideals in $R_{\mathbf{r},0}^+$.}
$ $

Throughout this section, we fix $t \in \gc 0,n \dc$ and $\mathbf{r}=(r_1,\dots,r_t) \in \mathbf{R}_t$, and we define 
$w_{\mathbf{r}}$ as in the previous section.

As in \cite[2.1]{glen1}, we set $\displaystyle{R_{\mathbf{r},0}^+=\frac{R}{\langle \yia \mbox{ $\mid$ } (\ia)
 \in w_{\mathbf{r}} \rangle}}$.
\\$ $

Recall (See \cite{glen1}, 2.1) that $R_{\mathbf{r},0}^+ $ can be written as an iterated Ore extension over $\mathbb{K}$. 
Thus, $R_{\mathbf{r},0}^+ $ is a Noetherian domain. Moreover, since $q$ is not a root of unity, 
it follows from \cite[Theorem 3.2]{gl4} that all primes 
of $R$ are completely prime and thus, since this property survives in factors, all primes in the algebra 
$R_{\mathbf{r},0}^+$ are completely prime. 

 Observe now that, since the indeterminates $\yia$ are $\hc$-eigenvectors, 
$\langle \yia \mbox{ $\mid$ } (\ia) \in w_{\mathbf{r}} \rangle$ 
is an $\hc$-invariant ideal of $R$. Hence, the action of $\hc$ on $R$ induces an action of $\hc$ 
on $R_{\mathbf{r},0}^+$ by automorphisms.  As usually, an \textbf{$\hc $-eigenvector} $x$ of $R_{\mathbf{r},0}^+$ 
is a nonzero element $x \in R_{\mathbf{r},0}^+$ such that $h(x) \in \mathbb{K}^*x$ 
for each $h \in \hc$, and an ideal $I$ of $R_{\mathbf{r},0}^+$ is said to be \textbf{$\hc$-invariant}
 if $h(I) =I$ for all $h\in \hc$. Further, we denote by \underline{$\hc$-$Spec(R_{\mathbf{r},0}^+)$} 
the set of $\hc$-invariant prime ideals of $R_{\mathbf{r},0}^+$.
\\$ $

\begin{notas}
$ $
\begin{itemize}
\item We denote by $\pi_{\mathbf{r},0}^+ : R \rightarrow R_{\mathbf{r},0}^+$ the canonical 
surjective $\mathbb{K}$-algebra homomorphism.
\item If $(\ia) \in \gc 1,n \dc^2$, $\overline{Y}_{\ia}$ denotes the element of $R_{\mathbf{r},0}^+$ 
defined by $\overline{Y}_{\ia}:=\pi_{\mathbf{r},0}^+ (\yia)$.
\end{itemize}
\end{notas}
$ $

Let $(\gamma_1,\dots,\gamma_n) \in \Gamma_{\mathbf{r}}$ (See Notation \ref{notagammar}) and 
define $w_{\mathbf{r},(\gamma_1,\dots,\gamma_n)}$ as in Theorem \ref{theofond}. 
Recall (See Theorem \ref{theofond}) that $w_{\mathbf{r},(\gamma_1,\dots,\gamma_n)}$ is an element of $W$ and that the 
$\hc$-invariant prime ideal $J_{w_{\mathbf{r},(\gamma_1,\dots,\gamma_n)}}$ of $R$ contains the indeterminates 
$\yia$ with $ (\ia) \in w_{\mathbf{r}}$, so that $\langle \yia \mbox{ $\mid$ } (\ia) \in w_{\mathbf{r}} \rangle 
\subseteq J_{w_{\mathbf{r},(\gamma_1,\dots,\gamma_n)}}$. Thus, 
$\pi_{\mathbf{r},0}^+ \left( J_{w_{\mathbf{r},(\gamma_1,\dots,\gamma_n)}} \right)$ is a (completely) prime ideal 
of $R_{\mathbf{r},0}^+$. More precisely, we have:
\\$ $

\begin{prop} \label{theofond2}
$ $
\\$J_{w_{\mathbf{r},(\gamma_1,\dots,\gamma_n)}}^+ := 
\pi_{\mathbf{r},0}^+ \left( J_{w_{\mathbf{r},(\gamma_1,\dots,\gamma_n)}} \right)$ is an $\hc$-invariant (completely) 
prime ideal of $R_{\mathbf{r},0}^+$ which does not contain the $\overline{Y}_{r_k,k}$ ($k \in \gc 1,t \dc$).
\end{prop} 
\preuve We have already explained that $J_{w_{\mathbf{r},(\gamma_1,\dots,\gamma_n)}}^+$ 
is a (completely) prime ideal of $R_{\mathbf{r},0}^+$. Moreover, since 
$J_{w_{\mathbf{r},(\gamma_1,\dots,\gamma_n)}}$ is $\hc$-invariant, it is easy to check that 
$J_{w_{\mathbf{r},(\gamma_1,\dots,\gamma_n)}}^+$ is also $\hc$-invariant. Finally, 
since $J_{w_{\mathbf{r},(\gamma_1,\dots,\gamma_n)}}$ does not contain the indeterminates $Y_{r_k,k}$ with 
$k \in \gc 1,t \dc$ (See Theorem \ref{theofond}), $J_{w_{\mathbf{r},(\gamma_1,\dots,\gamma_n)}}^+$ does not 
contain the $\overline{Y}_{r_k,k}=\pi_{\mathbf{r},0}^+ (Y_{r_k,k})$ with $k\in \gc 1,t \dc$. \fin
\\$ $

\subsection{$\hc$-invariant prime ideals in $R_{\mathbf{r}}^+$.}
\label{sectionR+r}
$ $

As in the previous section, we fix $t \in \gc 0,n \dc$ and $\mathbf{r}=(r_1,\dots,r_t) \in \mathbf{R}_t$. 
In \cite[2.1]{glen1}, Goodearl and Lenagan have observed that the $\overline{Y}_{r_k,k}$ with $k \in \gc 1,t \dc$ 
are regular normal elements in $R_{\mathbf{r},0}^+$, so that we can form the Ore localization: 
$$R_{\mathbf{r}}^+:=R_{\mathbf{r},0}^+ S_{\mathbf{r}}^{-1},$$
where $S_{\mathbf{r}}$ denotes the multiplicative system of $R_{\mathbf{r},0}^+$ 
generated by the $\overline{Y}_{r_k,k}$ with $k \in \gc 1,t \dc$. 
\\$ $

In the previous section, we have noted that all the primes of $R_{\mathbf{r},0}^+$ are completely prime. 
Since this property survives in localization, all the primes of $R_{\mathbf{r}}^+$ are 
also completely prime.

Observe now that, since the $\overline{Y}_{r_k,k}$ with $k\in \gc 1,t \dc$ are $\hc$-eigenvectors of 
$R_{\mathbf{r},0}^+$, the action of $\hc$ on $R_{\mathbf{r},0}^+$ extends to an action 
of $\hc$ on $R_{\mathbf{r}}^+$ by automorphisms. We say that an ideal $I$ of $R_{\mathbf{r}}^+$ is 
\textbf{$\hc$-invariant} if $h(I) =I$ for all $h\in \hc$ 
and we denote by \underline{$\hc$-$Spec(R_{\mathbf{r}}^+)$} the set of $\hc$-invariant prime ideals of $R_{\mathbf{r}}^+$. 
Observe now that contraction and extension provide inverse bijections 
between the set $\hc$-$Spec(R_{\mathbf{r}}^+)$ and the set of those 
$\hc$-invariant prime ideals of $R_{\mathbf{r},0}^+$ which are disjoint from $S_{\mathbf{r}}$.

Let $(\gamma_1,\dots,\gamma_n) \in \Gamma_{\mathbf{r}}$ (See Notation \ref{notagammar}) and define 
$\displaystyle{w_{\mathbf{r},(\gamma_1,\dots,\gamma_n)}}$ as in Theorem \ref{theofond}. 
By Proposition \ref{theofond2}, $J_{w_{\mathbf{r},(\gamma_1,\dots,\gamma_n)}}^+ := 
\pi_{\mathbf{r},0}^+ \left( J_{w_{\mathbf{r},(\gamma_1,\dots,\gamma_n)}} \right)$ is an $\hc$-invariant (completely) 
prime ideal of $R_{\mathbf{r},0}^+$ which does not contain the $\overline{Y}_{r_k,k}$ ($k \in \gc 1,t \dc$). Since 
$S_{\mathbf{r}}$ is generated by the $\overline{Y}_{r_k,k}$ ($k \in \gc 1,t \dc$), 
$J_{w_{\mathbf{r},(\gamma_1,\dots,\gamma_n)}}^+ $ is an $\hc$-invariant (completely) 
prime ideal of $R_{\mathbf{r},0}^+$ which is disjoint from $S_{\mathbf{r}}$. Thus, we have the 
following statement:
\\$ $

\begin{prop}
\label{theofond3}
$ $
\\$J_{w_{\mathbf{r},(\gamma_1,\dots,\gamma_n)}}^+ S_{\mathbf{r}}^{-1}$ is an $\hc$-invariant (completely) 
prime ideal of $R_{\mathbf{r}}^+$.
\end{prop} 
$ $

We will prove later (See Section \ref{lastsection}) that the $J_{w_{\mathbf{r},(\gamma_1,\dots,\gamma_n)}}^+ S_{\mathbf{r}}^{-1}$ 
($(\gamma_1,\dots,\gamma_n) \in \Gamma_{\mathbf{r}}$) are exactly the $\hc$-invariant prime ideals of $R_{\mathbf{r}}^+$.
\\$ $

We deduce from the above Proposition \ref{theofond3} that:\\
$ $

\begin{cor}
\label{corfond3}
$ $
\\$R_{\mathbf{r}}^+$ has at least $1^{r_1}2^{r_2-r_1}\dots t^{r_t-r_{t-1}} (t+1)^{n-r_t}$ 
$\hc$-invariant prime ideals.
\end{cor} 
\preuve It follows from Proposition \ref{theofond3} that 
$R_{\mathbf{r}}^+$ has at least $\mid \Gamma_{\mathbf{r}} \mid$ $\hc$-invariant prime ideals, and 
it is obvious that $\mid \Gamma_{\mathbf{r}} \mid \ = 1^{r_1}2^{r_2-r_1}\dots t^{r_t-r_{t-1}} (t+1)^{n-r_t}. \mbox{ \fin}$ 
\\$ $

\subsection{Number of rank $t$ $\hc$-invariant prime ideals in $\mnk$.}
\label{sectioncount}
$ $

For convenience, we recall the following definitions (See \cite{pw}):
\\$ $
 \begin{defis}
$ $
\begin{itemize}
\item Let $m$ be a positive integer and let $M=(x_{\ia})_{(\ia) \in \gc 1,m \dc ^2}$ be a square $q$-quantum matrix.
\\The quantum determinant of $M$ is defined by:
$$det_q(M):=\sum_{\sigma \in S_m}(-q)^{l(\sigma)}x_{1,\sigma(1)} \dots
x_{m,\sigma(m)},$$
where $S_m$ denotes the group of permutations of $\gc 1,m \dc$ and 
$l(\sigma)$ denotes the length of the $m$-permuation $\sigma$.
\item Let $\mathcal{Y}:=(Y_{\ia})_{(\ia) \in \gc 1,n \dc^2 }$ be the $q$-quantum matrix of the canonical 
generators of $R$.
\\The quantum determinant of a square sub-matrix of $\mathcal{Y}$ is called a quantum minor. 
 \end{itemize}
\end{defis}
$ $

We can now define the rank $t$ $\hc$-invariant prime ideals of $R$, as follows:
\\$ $

\begin{defi}
$ $
\\Let $t \in \gc 0,n \dc$. An $\hc$-invariant prime ideal $J$ of $R=\mnk$ has rank $t$ 
if $J$ contains all $(t+1) \times (t+1)$ quantum minors but not all $t \times t$ quantum minors.
\\As in \cite[3.6]{glen1}, we denote by $\hc$-$Spec^{[t]} (R)$ 
the set of rank $t$ $\hc$-invariant prime ideals of $R$.
\end{defi}
$ $

Note that there is only one element in $\hc$-$Spec^{[0]} (R)$: $\langle \yia \mbox{ $\mid$ } 
(\ia) \in \gc 1,n \dc^2 \rangle$, the augmentation ideal of $R$. Further, 
Goodearl and Lenagan have observed (See \cite{glen1}, 3.6) that 
$\mid \hc \mbox{-}Spec^{[1]} (R) \mid \ =(2^n -1)^2$ and $\mid \hc \mbox{-}Spec^{[n]} (R) \mid \ =(n!)^2$.
\\$ $

\begin{obs}
\label{obspartition}
$ $
\\The sets $\hc$-$Spec^{[t]} (R)$ ($t \in \gc 0,n \dc$) partition the set $\hc$-$Spec^{[t]} (R)$ .
\end{obs} 
\preuve Let $P$ be an $\hc$-invariant prime ideal of $R$. Let $t \in \gc 0,n \dc$ be maximal such that $P$ 
does not contain all $t \times t$ quantum minors. Then $P$ clearly belongs 
to $\hc$-$Spec^{[t]} (R)$. Hence, we have proved that 
$\displaystyle{\hc \mbox{-} Spec (R) = \bigcup_{t \in \gc 0,n \dc} \hc \mbox{-}Spec^{[t]} (R)}$. Since this union is obviously 
disjoint, we get $\displaystyle{\hc \mbox{-}Spec (R) = \bigsqcup_{t \in \gc 0,n \dc} \hc \mbox{-}Spec^{[t]} (R)}$, 
as desired. \fin
\\$ $

In \cite{glen1}, the authors have established the following result 
that will be our starting point to compute the cardinality of $\hc$-$Spec^{[t]} (R)$:
\\$ $

\begin{prop}[See \cite{glen1}, 3.6]
$ $
\\For all $t \in \gc 0,n \dc$, we have 
$\displaystyle{\mid \hc\mbox{-}Spec^{[t]} (R) \mid \ = 
\left( \sum_{\mathbf{r} \in \mathbf{R}_t} \mid \hc \mbox{-} Spec(R_{\mathbf{r}}^+) \mid \right)^2}$.
\end{prop}
$ $
 
Before computing $\mid \hc\mbox{-}Spec^{[t]} (R) \mid$, we first give a lower bound for 
$ \displaystyle{\sum_{\mathbf{r} \in \mathbf{R}_t} \mid \hc \mbox{-} Spec(R_{\mathbf{r}}^+) \mid }$.
\\$ $

\begin{prop}
\label{propmino}
$ $
\\For any $t \in \gc 0,n \dc$, we have 
$$\sum_{\mathbf{r} \in \mathbf{R}_t} \mid \hc \mbox{-} Spec(R_{\mathbf{r}}^+) \mid \ \geq t! S(n+1,t+1),$$ 
where $S(n+1,t+1)$ denotes the Stirling number of second kind associated to $n+1$ and $t+1$ 
(See, for instance, \cite{stanleybook} for the definition of $S(n+1,t+1)$).
\end{prop} 
\preuve First, we deduce from Corollary \ref{corfond3} the following inequality: 
\begin{eqnarray}
\label{equationdebut}
\sum_{\mathbf{r} \in \mathbf{R}_t} \mid \hc \mbox{-} Spec(R_{\mathbf{r}}^+) \mid \ \geq  
\sum_{\mathbf{r} \in \mathbf{R}_t} 1^{r_1}2^{r_2-r_1}\dots t^{r_t-r_{t-1}} (t+1)^{n-r_t}. 
\end{eqnarray}
On the other hand, we know (See \cite{stanleybook}, Exercise 16 p46) that: 
\begin{eqnarray}
\label{equationutile1}
S(n+1,t+1) = \sum_{a_1+\dots+a_{t+1}=n+1}1^{a_1-1} 2^{a_2-1} \dots (t+1)^{a_{t+1}-1}.
\end{eqnarray}
Observe now that the map $f:\{(a_1,\dots,a_{t+1})\in (\mathbb{N}^*)^{t+1} \mbox{ $\mid$ }a_1+\dots+a_{t+1}=n+1\} 
\rightarrow \{(r_1,\dots,r_{t})\in (\mathbb{N}^*)^{t} \mbox{ $\mid$ }1 \leq r_1 < \dots < r_{t} \leq n\}=
 \mathbf{R}_t$ defined by $f(a_1,\dots,a_{t+1})=(a_1,a_1+a_2,\dots,a_1+\dots+a_{t})$ is a bijection and that 
its inverse $f^{-1}$ is defined by $f^{-1}(r_1,\dots,r_t)=(r_1,r_2-r_1,\dots,r_t-r_{t-1},n+1-r_t)$ for all 
$(r_1,\dots,r_t) \in \mathbf{R}_t$. Thus, by means of the change of variables 
$(a_1,\dots,a_{t+1})=f^{-1}(r_1,\dots,r_t)$, the above equality (\ref{equationutile1}) is 
transformed to 
$$S(n+1,t+1) = \sum_{1 \leq r_1 < \dots < r_{t} \leq n}1^{r_1-1} 2^{r_2-r_1-1} \dots 
t^{r_t-r_{t-1}-1} (t+1)^{n-r_{t}},$$
so that 
$$t! S(n+1,t+1)= \sum_{(r_1,\dots,r_{t}) \in \mathbf{R}_t}1^{r_1} 2^{r_2-r_1} \dots 
t^{r_t-r_{t-1}} (t+1)^{n-r_{t}}.$$ 
Thus, we deduce from inequality (\ref{equationdebut}) that:
$$\sum_{\mathbf{r} \in \mathbf{R}_t} \mid \hc \mbox{-} Spec(R_{\mathbf{r}}^+) \mid \ \geq  
t!S(n+1,t+1), $$
as desired. \fin
\\$ $

\begin{rem}
\label{remmino}
$ $
\\The proof of the above Proposition \ref{propmino} shows that, if there exists 
$t \in \gc 0,n \dc$ and $\mathbf{r}=(r_1,\dots,r_t) \in \mathbf{R}_t$ such that 
$\mid \hc \mbox{-} Spec(R_{\mathbf{r}}^+) \mid \ > 1^{r_1}2^{r_2-r_1}\dots t^{r_t-r_{t-1}} (t+1)^{n-r_t}$, then 
 $$\sum_{\mathbf{r} \in \mathbf{R}_t} \mid \hc \mbox{-} Spec(R_{\mathbf{r}}^+) \mid \ > t! S(n+1,t+1).$$
\end{rem}
$ $

We can now prove our main result which was conjectured by Goodearl, Lenagan and McCammond:
\\$ $

\begin{theo}
\label{numberrankt}
$ $
\\If $t \in \gc 0,n \dc$, then $\mid \hc\mbox{-}Spec^{[t]} (R) \mid \ = 
\left( t! S(n+1,t+1) \right)^2$.
\end{theo}
\preuve First, since the sets $\hc\mbox{-}Spec^{[t]} (R)$ ($t \in \gc 0,n \dc$) partition $\hc$-$Spec(R)$ 
(See Observation \ref{obspartition}), we have :
$$ \mid \hc \mbox{-} Spec(R) \mid \ = \sum_{t=0}^n \mid \hc\mbox{-}Spec^{[t]} (R) \mid.$$
Recall now (See Proposition \ref{propnbreHprem}) that $\mid \hc \mbox{-} Spec(R) \mid$ is equal to the poly-Bernoulli
 number $B_n^{(-n)}$. Thus, we deduce from the above equality that:
$$ B_n^{(-n)} = \sum_{t=0}^n \mid \hc\mbox{-}Spec^{[t]} (R) \mid.$$
Further, by \cite[Theorem 2]{kaneko2}, $B_n^{(-n)}$ can also be written as follows:
$$B_n^{(-n)}=  \sum_{t=0}^n \left( t! S(n+1,t+1) \right)^2.$$ Hence, we have:
$$\sum_{t=0}^n \mid \hc\mbox{-}Spec^{[t]} (R) \mid \ =
\sum_{t=0}^n \left( t! S(n+1,t+1) \right)^2,$$ 
that is:
\begin{eqnarray}
\label{equationpolyBer}
\sum_{t=0}^n  \left( \mid \hc\mbox{-}Spec^{[t]} (R) \mid \ - \left( t! S(n+1,t+1) \right)^2 \right)=0. 
\end{eqnarray}

On the other hand, recall (See \cite{glen1}, 3.6) that $ \displaystyle{\mid \hc\mbox{-}Spec^{[t]} (R) \mid \ = 
\left( \sum_{\mathbf{r} \in \mathbf{R}_t} \mid \hc \mbox{-} Spec(R_{\mathbf{r}}^+) \mid \right)^2}$. Thus, 
since $\displaystyle{\sum_{\mathbf{r} \in \mathbf{R}_t} \mid \hc \mbox{-} Spec(R_{\mathbf{r}}^+) \mid \ \geq t! S(n+1,t+1)}$ 
(See Proposition \ref{propmino}), we have: 
$$\mid \hc\mbox{-}Spec^{[t]} (R) \mid \ 
\geq  \left( t! S(n+1,t+1) \right)^2.$$
In other words, each of the terms which appears in the sum on the left hand side of (\ref{equationpolyBer}) is 
non-negative. Since this sum is equal to zero, each term of this sum must be zero, that is, 
 for all $t \in \gc 0,n \dc$, we have:
$$\mid \hc\mbox{-}Spec^{[t]} (R) \mid \ = 
\left( t! S(n+1,t+1) \right)^2. \mbox{ \fin}$$
$ $

\begin{rem}
$ $
\\The cases $t=0$, $t=1$ and $t=n$ were already known (See \cite{glen1}, 3.6).
\end{rem}
$ $

\subsection{Description of the set $\hc$-$Spec(R_{\mathbf{r}}^+)$.}
\label{lastsection}
$ $

Throughout this section, we fix $t \in \gc 0,n \dc$ and $\mathbf{r}=(r_1,\dots,r_t) \in \mathbf{R}_t$. 
We now use the above Theorem \ref{numberrankt} to obtain a description of the set $\hc$-$Spec(R_{\mathbf{r}}^+)$. 
More precisely, we show that the only $\hc$-invariant prime ideals of $R_{\mathbf{r}}^+$ are 
those obtained in Proposition \ref{theofond3}, that is, in the notations of Section \ref{sectionR+r}:
\\$ $

\begin{theo}
$$ \hc \mbox{-}Spec(R_{\mathbf{r}}^+)=\{ J_{w_{\mathbf{r},(\gamma_1,\dots,\gamma_n)}}^+ S_{\mathbf{r}}^{-1} 
\mbox{ $\mid$ } (\gamma_1,\dots,\gamma_n) \in \Gamma_{\mathbf{r}} \}.$$
\end{theo} 
\preuve We already know (See Proposition \ref{theofond3}) that 
$$ \hc \mbox{-}Spec(R_{\mathbf{r}}^+) \supseteq \{ J_{w_{\mathbf{r},(\gamma_1,\dots,\gamma_n)}}^+ S_{\mathbf{r}}^{-1} 
\mbox{ $\mid$ } (\gamma_1,\dots,\gamma_n) \in \Gamma_{\mathbf{r}} \}.$$
Assume now that $$ \hc \mbox{-}Spec(R_{\mathbf{r}}^+) \varsupsetneqq \{ J_{w_{\mathbf{r},(\gamma_1,\dots,\gamma_n)}}^+ S_{\mathbf{r}}^{-1} 
\mbox{ $\mid$ } (\gamma_1,\dots,\gamma_n) \in \Gamma_{\mathbf{r}} \}.$$
Then we have $\mid  \hc \mbox{-}Spec(R_{\mathbf{r}}^+) \mid \ > \ \mid \Gamma_{\mathbf{r}} \mid$. Since 
$\mid \Gamma_{\mathbf{r}} \mid \ = 1^{r_1}2^{r_2-r_1}\dots t^{r_t-r_{t-1}} (t+1)^{n-r_t} $, we get 
$\mid \hc \mbox{-}Spec(R_{\mathbf{r}}^+) \mid \ >  1^{r_1}2^{r_2-r_1}\dots t^{r_t-r_{t-1}} (t+1)^{n-r_t} $. 
Thus, it follows from Remark \ref{remmino} that 
$$\sum_{\mathbf{r} \in \mathbf{R}_t} \mid \hc \mbox{-} Spec(R_{\mathbf{r}}^+) \mid \ > t! S(n+1,t+1).$$
Hence we have 
$$\left( \sum_{\mathbf{r} \in \mathbf{R}_t} \mid \hc \mbox{-} Spec(R_{\mathbf{r}}^+) \mid \right)^2 > 
\left( t! S(n+1,t+1) \right)^2.$$
Recall now (See \cite[3.6]{glen1}) that 
$$\mid \hc\mbox{-}Spec^{[t]} (R) \mid \ = 
\left( \sum_{\mathbf{r} \in \mathbf{R}_t} \mid \hc \mbox{-} Spec(R_{\mathbf{r}}^+) \mid \right)^2.$$
All this together shows that $\mid \hc\mbox{-}Spec^{[t]} (R) \mid \ > \left( t! S(n+1,t+1) \right)^2$.

However, it follows from Theorem \ref{numberrankt} that 
$\mid \hc\mbox{-}Spec^{[t]} (R) \mid \ = \left( t! S(n+1,t+1) \right)^2$. 
This is a contradiction and thus we have proved that 
$ \hc \mbox{-}Spec(R_{\mathbf{r}}^+)=\{ J_{w_{\mathbf{r},(\gamma_1,\dots,\gamma_n)}}^+ S_{\mathbf{r}}^{-1} 
\mbox{ $\mid$ } (\gamma_1,\dots,\gamma_n) \in \Gamma_{\mathbf{r}} \}$. \fin
\\$ $

\begin{flushleft}
\textbf{Acknowledgments.}
\end{flushleft}
I thank T.H. Lenagan for very helpful conversations, and K.R. Goodearl for useful comments.
\\$ $

\def\refname{References}

\bibliographystyle{amsplain}
\bibliography{biblio}

\providecommand{\bysame}{\leavevmode\hbox to3em{\hrulefill}\thinspace}
\providecommand{\MR}{\relax\ifhmode\unskip\space\fi MR }
\providecommand{\MRhref}[2]{%
  \href{http://www.ams.org/mathscinet-getitem?mr=#1}{#2}
}
\providecommand{\href}[2]{#2}
\begin{thebibliography}{10}

\bibitem{kaneko2}
T.~Arakawa and M.~Kaneko, \emph{On poly-{B}ernoulli numbers}, Comment Math.
  Univ. St. Paul \textbf{48 (2)} (1999), 159--167.

\bibitem{c3}
G.~Cauchon, \emph{Quotients premiers de ${O}_q \left( \mathcal{M}_n (k)
  \right)$}, J. Algebra \textbf{180} (1996), 530--545.

\bibitem{c1}
\bysame, \emph{Effacement des d\'erivations et spectres premiers d'alg\`ebres
  quantiques}, J. Algebra. \textbf{260} (2003), 476--518.

\bibitem{c2}
\bysame, \emph{Spectre premier de ${O}_q \left( \mathcal{M}_n (k) \right)$,
  image canonique et s\'eparation normale}, J. Algebra. \textbf{260} (2003),
  519--569.

\bibitem{glen1}
K.R. Goodearl and T.H. Lenagan, \emph{Prime ideals invariant under winding
  automorphisms in quantum matrices}, Internat. J. Math. \textbf{13} (2002),
  497--532.

\bibitem{glen2}
\bysame, \emph{Winding-invariant prime ideals in quantum $3 \times 3$
  matrices}, J. Algebra. \textbf{260} (2003), 657--687.

\bibitem{gl4}
K.R. Goodearl and E.S. Letzter, \emph{Prime factor algebras of the coordinate
  ring of quantum matrices}, Proc. Amer. Math. Soc. \textbf{121} (1994),
  1017--1025.

\bibitem{gl2}
\bysame, \emph{Prime and primitive spectra of multiparameter quantum affine
  spaces, \textnormal{Trends in ring theory (Miskolc, 1996)}}, Canad. Math.
  Soc. Conf. Proc. Series, vol.~22, 1998, pp.~39--58.

\bibitem{gl1}
\bysame, \emph{The {D}ixmier-{M}oeglin equivalence in quantum coordinate rings
  and quantized {W}eyl algebras}, Trans. Amer. Math. Soc. \textbf{352} (2000),
  1381--1403.

\bibitem{kaneko1}
M.~Kaneko, \emph{Poly-{B}ernoulli numbers}, J. de Th\'eorie des Nombres de
  Bordeaux \textbf{9} (1997), 221--228.

\bibitem{lau2}
S.~Launois, \emph{Generators for the $\hc$-invariant prime ideals in $\mmpc$},
  to appear in Proceedings of the Edinburgh Mathematical Society.

\bibitem{lau}
\bysame, \emph{Les id\'eaux premiers invariants de $\mmpc$}, to appear in J.
  Algebra.

\bibitem{lau3}
\bysame, \emph{Id\'eaux premiers $\hc$-invariants de l'alg\`ebre des matrices
  quantiques}, Th\`ese de doctorat, Universit\'e de Reims, 2003.

\bibitem{pw}
B.J. Parshall and J.P. Wang, \emph{Quantum {L}inear {G}roups}, Mem. Amer. Math.
  Soc., 439, 1991.

\bibitem{stanleybook}
R.P. Stanley, \emph{Enumerative {C}ombinatorics {I}}, Cambridge University
  Press, 1997.

\end{thebibliography}

\end{document}